\documentclass{article}
\usepackage{amsmath}
\usepackage{amssymb}
\usepackage[a4paper, margin={1.2in}]{geometry}
\usepackage[utf8]{inputenc}
\usepackage{lipsum}

\newcommand\blfootnote[1]{%
  \begingroup
  \renewcommand\thefootnote{}\footnote{#1}%
  \addtocounter{footnote}{-1}%
  \endgroup
}

\begin{document}

\title{Indiscernibles and satisfaction classes in arithmetic}
\author{Ali Enayat}
\maketitle

\begin{abstract}


\noindent We investigate the theory $\mathrm{PAI}$ (Peano Arithmetic with
Indiscernibles). Models of $\mathrm{PAI}$\ are of the form $(\mathcal{M},I)$%
, where $\mathcal{M}$\ is a model of $\mathrm{PA}$, $I$ is an unbounded set
of order indiscernibles over $\mathcal{M}$, and $(\mathcal{M},I)$\ satisfies
the extended induction scheme for formulae mentioning $I$. Our main results
are Theorems A and B below.\medskip

\noindent \textbf{Theorem A}.~\textit{Let} $\mathcal{M}$\ \textit{be a
nonstandard model of} $\mathrm{PA}$\textit{\ of any cardinality}. $\mathcal{M%
}$ \textit{has an expansion to a model of }$\mathrm{PAI}$\ \textit{iff} $%
\mathcal{M}$ \textit{has an inductive partial satisfaction class.}\medskip

\noindent Theorem A yields the following corollary, which provides a new
characterization of countable recursively saturated models of $\mathrm{PA}$:
\medskip

\noindent \textbf{Corollary.}~\textit{A countable model} $\mathcal{M}$
\textit{of} $\mathrm{PA}$ \textit{is recursively saturated iff }$\mathcal{M}$
\textit{has an expansion to a model of }$\mathrm{PAI}$.\medskip

\noindent \textbf{Theorem B}.~\textit{There is a sentence }$\alpha $ \textit{%
in the language obtained by adding a unary predicate} $I(x)$ \textit{to the
language of arithmetic such that given any nonstandard model }$\mathcal{M}$\
\textit{of} $\mathrm{PA}$\textit{\ of any cardinality}, $\mathcal{M}$
\textit{has an expansion to a model of }$\mathrm{PAI}+\alpha $\ \textit{iff}
$\mathcal{M}$ \textit{has a inductive full satisfaction class.}\medskip
\end{abstract}

\begin{center}
\bigskip

\textbf{1.~INTRODUCTION\bigskip }
\end{center}

We investigate an extension of $\mathrm{PA}$ (Peano Arithmetic), denoted $%
\mathrm{PAI}$, which is equipped with a designated unbounded class of
indiscernibles (see Section 3 for the precise definition). The motivation to
study $\mathrm{PAI}$ arose from the study \cite{Ali-ZFI} of the
set-theoretic counterpart $\mathrm{ZFI}_{<}$ of $\mathrm{PAI,}$ where it is
shown that there is an intimate relationship between $\mathrm{ZFI}_{<}$ and
large cardinals, thus indicating that the set-theoretical consequences of $%
\mathrm{ZFI}_{<}$ go well beyond $\mathrm{ZFC}$. \medskip

In light of the results obtained in \cite{Ali-ZFI} it is natural to
investigate \textrm{PAI} since it is well-known \cite{Richard-Tin Lok} that $%
\mathrm{PA}$ is bi-interpretable with the theory \textrm{ZF}$^{-\infty }+%
\mathrm{TC,}$ where \textrm{ZF}$^{-\infty }$ is the system of set theory
obtained from \textrm{ZF} by replacing the axiom of infinity by its
negation, and $\mathrm{TC}$ is the sentence asserting that every set is
contained in a transitive set (which in the presence of the other axioms
implies that the transitive closure of every set exists). The aforementioned
proof of the bi-interpretability of $\mathrm{PA}$ and \textrm{ZF}$^{-\infty
}+\mathrm{TC}$ can be readily extended to show the bi-interpretability of $%
\mathrm{PAI}$ and \textrm{ZFI}$^{-\infty }+\mathrm{TC.}$\medskip

Our main results are Theorems A and B of the abstract that relate $\mathrm{%
PAI}$ to the well-studied notions of (a) inductive partial satisfaction
classes and (b) inductive full satisfaction classes, which are intimately
connected (respectively) with the axiomatic theories of truth known as $%
\mathrm{UTB}$ and $\mathrm{CT}$ (see Section 2.2). After presenting
preliminaries in Section 2, (refinements of) Theorems A and B are
established in Section 3. In Section 4 we examine $\mathrm{PAI}$ through the
lens of interpretability, and in Section 5 we probe the model-theoretic
behavior of fragments of $\mathrm{PAI.}$\footnote{%
I am grateful to Lawrence Wong, Bartosz Wcis\l o, Mateusz \L e\l yk, Jim
Schmerl, Roman Kossak, Kentaro Fujimoto, Cezary Cie\'{s}li\'{n}ski, Lev
Beklemishev, and Athar Abdul-Quader (in reverse alphabetical order) for
their feedback and keen interest in this work. The research presented in
this paper was partially supported by the National Science Centre, Poland
(NCN), grant number 2019/34/A/HS1/00399.}\bigskip

\begin{center}
\textbf{2.~PRELIMINARIES}\bigskip
\end{center}

In this section we present the relevant notations, conventions, definitions,
and results that are needed in the subsequent sections.\textbf{\medskip }

\begin{center}
\noindent \textbf{2.1.}~\textbf{Theories and models\medskip }
\end{center}

\noindent \textbf{2.1.1.}~\textbf{Definition.}~The language of arithmetic, $%
\mathcal{L}_{\mathrm{A}}$, is $\{+,\cdot ,\mathrm{S},<,0\}.$ We use the
convention of writing $M$, $M_{0},$ $N$, etc. to (respectively) denote the
universes of discourse of structures $\mathcal{M}$, $\mathcal{M}_{0},$ $%
\mathcal{N},$ etc. In (a) through (g) below, $\mathcal{L}\supseteq \mathcal{L%
}_{\mathrm{A}}$ and $\mathcal{M}$ and $\mathcal{N}$ are $\mathcal{L}$%
-structures.\textbf{\medskip }

\noindent \textbf{(a)} $\Sigma _{0}(\mathcal{L})=\Pi _{0}=\Delta _{0}(%
\mathcal{L})$ = the collection of $\mathcal{L}$-formulae all of whose
quantifiers are bounded by $\mathcal{L}$-terms (i.e., they are of the form $%
\exists x\leq t,$ or of the form $\forall x\leq t,$ where $t$ is an $%
\mathcal{L}$-term not involving $x$. More generally, $\Sigma _{n+1}(\mathcal{%
L)}$ consists of formulae of the form $\exists x_{0}\cdot \cdot \cdot
\exists x_{k-1}\ \varphi $, where $\varphi \in \Pi _{n}(\mathcal{L)};$ and $%
\Pi _{n+1}(\mathcal{L)}$ consists of formulae of the form $\forall
x_{0}\cdot \cdot \cdot \forall x_{k-1}\ \varphi $, where $\varphi \in \Sigma
_{n}$ (with the convention that $k=0$ corresponds to an empty block of
quantifiers). \textit{We shall omit the reference to }$\mathcal{L}$\textit{\
if }$\mathcal{L}=\mathcal{L}_{\mathrm{A}}$\textit{, e.g., }$\Sigma
_{n}:=\Sigma _{n}(\mathcal{L}_{\mathrm{A}}\mathcal{)}$. Also, we shall write
$\Sigma _{n}(X\mathcal{)}$ instead of $\Sigma _{n}(\mathcal{L}_{\mathrm{A}}%
\mathcal{\cup \{}X\mathcal{\})},$ where $X$ is a new predicate symbol. We
often conflate formal symbols with their denotations (if there is no risk of
confusion).\textbf{\medskip }

\noindent \textbf{(b)} $\mathrm{PA}$ (Peano Arithmetic)\ is the result of
adding the scheme of induction for all $\mathcal{L}_{\mathrm{A}}$-formulae
to the finitely axiomatizable theory known as (Robinson's) $\mathrm{Q}$. $%
\mathrm{PA}(\mathcal{L})$ is the theory obtained by augmenting $\mathrm{PA}$
with the scheme of induction for all $\mathcal{L}$-formulae. $\mathrm{I}%
\Sigma _{n}(\mathcal{L})$ is the fragment of $\mathrm{PA}(\mathcal{L})$ with
the induction scheme limited to $\Sigma _{n}(\mathcal{L})$-formulae. Given a
new predicate $X$, we write $\mathrm{PA}(X)$\ and $\mathrm{I}\Sigma _{n}(X)$
(respectively) instead of $\mathrm{PA}(\mathcal{L}_{\mathrm{A}}\mathcal{\cup
\{}X\mathcal{\}})$ and $\mathrm{I}\Sigma _{n}(\mathcal{L}_{\mathrm{A}}%
\mathcal{\cup \{}X\mathcal{\}}).$\textbf{\medskip }

\noindent \textbf{(c)} If $\varphi (x)$ is an $\mathcal{L}$-formula, $%
\varphi ^{\mathcal{M}}:=\{m\in M:\mathcal{M}\models $ $\varphi (x)\}$. For $%
X\subseteq M$, then we say that $X$\ is $\mathcal{M}$-\textit{definable} if $%
X$\ is first order definable (parameters allowed) in $\mathcal{M}$. \textbf{%
\medskip }

\noindent \textbf{(d)} A subset $X$ of $M$ is $\mathcal{M}$\textit{-finite }(%
\textit{or }$\mathcal{M}$\textit{-coded})\textit{\ }if $X=c_{E}$ for some $%
c\in M$, where $c_{E}=\{m\in M:\mathcal{M}\models m\in _{\mathrm{Ack}}c\},$
and $m\in _{\mathrm{Ack}}c$ is shorthand for the formula expressing
\textquotedblleft the $m$-th bit of the binary expansion of $c$ is
1\textquotedblright .\textbf{\medskip }

\noindent \textbf{(e) }A subset $X$\ of $M$ is said to be \textit{%
piecewise-coded} (in $\mathcal{M}$)\ if $\{x\in M:x<m$ and $x\in X\}$ is $%
\mathcal{M}$-finite for each $m\in M$.\textbf{\medskip }

\noindent \textbf{(f) }$\mathcal{M}$ is \textit{rather classless} if any
piecewise-coded subset of $M$ is already $\mathcal{M}$-definable (By a
theorem of Kaufmann, every extension of \textrm{PA} has a recursively
saturated rather classless model \cite[Theorem 10.1.5]{Roman-Jim's book}).%
\textbf{\medskip }

\noindent \textbf{(g)} We identify the longest well-founded initial submodel
of models of PA with the ordinal $\omega .$\textbf{\medskip }

The following result was established by Kossak \cite[Proposition 3.2]{Roman
(on certain)}\ for models $\mathcal{M}$ of $\mathrm{PA}$; the generalization
to models of $\mathrm{I}\Delta _{0}+\mathrm{Exp}$ appears in \cite[Lemma 4.2]%
{Ali-Fedor} (note that $\mathrm{Exp}$ is the axiom stating the totality of
the exponential function).\textbf{\medskip }

\noindent \textbf{2.1.2.~Theorem.}~\textit{Let }$\mathcal{M}\models \mathrm{I%
}\Delta _{0}+\mathrm{Exp},$ \textit{and} $X\subseteq M.$ \textit{The
following are equivalent}:\textbf{\medskip }

\noindent $(i)$ $(\mathcal{M},X)\models \mathrm{I}\Delta _{0}(X)$\textbf{%
.\medskip }

\noindent $(ii)$ $X$ \textit{is piecewise-coded in} $\mathcal{M}$.\textbf{%
\medskip }

\pagebreak

\begin{center}
\noindent \textbf{2.2.}~\textbf{Satisfaction classes, truth theories, and
recursive saturation\medskip }
\end{center}

\noindent \textbf{2.2.1.}~\textbf{Definition.}~Suppose $\mathcal{M}\models
\mathrm{PA}$, and $S\subseteq M.$ \textbf{\medskip }

\noindent \textbf{(a)} $S$ is said to be \textit{inductive}, if $(\mathcal{M}%
,S)\models \mathrm{PA}(S)$.\textbf{\medskip }

\noindent \textbf{(b)} $S$ is said to be a \textit{partial satisfaction class%
} if $S$ satisfies Tarski's recursive conditions for a satisfaction
predicate for all standard formulae. Thus a typical member of $S$ is of the
form $\left\langle \varphi ,\overline{a}\right\rangle ,$ where $\varphi \in
\mathrm{Form}_{m}^{\mathcal{M}}=$ the set of $\mathcal{L}_{\mathrm{A}}$%
-formulae in $\mathcal{M}$ with $m$ free variables, where $m\in M$ (note
that $\varphi $ need not be standard) and $\overline{a}\in M$ is an $m$%
-tuple in the sense of $\mathcal{M}$. \textbf{\medskip }

\noindent \textbf{(c)} $S$ is said to be a \textit{full satisfaction class}
if $S$ satisfies Tarski's recursive conditions for a satisfaction predicate
for all formulae in $\mathcal{M}$.

\begin{itemize}
\item For better readability we will often write $\left\langle \varphi ,%
\overline{a}\right\rangle \in S$ or $\varphi (\overline{a})\in S$ instead of
the more official $S(\left\langle \varphi ,\overline{a}\right\rangle ).$
Also, if $\varphi $ is a sentence (i.e., has no free variables), we will
write $\varphi \in S$ instead of $\left\langle \varphi ,\varnothing
\right\rangle \in S$ (where $\varnothing $ is the empty tuple).
\end{itemize}

The theories $\mathrm{UTB}$ (Uniform Tarski Biconditionals) and CT
(Compositional Truth) described below are well studied in the literature of
axiomatic theories of truth (see, e.g., the monographs by Cie\'{s}li\'{n}ski
\cite{Cezary book} and Halbach \cite{Halbach book}).

\begin{itemize}
\item Note that for the purposes of this paper, satisfaction and truth are
interchangeable, but in general there are subtle differences between the
two, see \cite{Cezary Tand S}.$\medskip $
\end{itemize}

\noindent \textbf{2.2.2.}~\textbf{Definition.}~In what follows $T(x)$ is a
new unary predicate$,$ $c$ is a new constant symbol, $\mathrm{Form}_{1}$ is
the set of (G\"{o}del-numbers of) $\mathcal{L}_{\mathrm{A}}$-formulae with
exactly one free variable, and $\overset{\cdot }{x}$ is the arithmetically
definable function that outputs the numeral for $x$ given the input $x$.$%
\medskip $

\noindent \textbf{(a)} $\mathrm{UTB}$ is $\mathrm{PA}(T)+\{\forall x(\varphi
(x)\leftrightarrow T(\ulcorner \varphi (\overset{\cdot }{x})\urcorner
):\varphi (x)\in \mathrm{Form}_{\mathrm{1}}\}$.$\medskip $

\noindent \textbf{(b)} $\mathrm{UTB}(c)=\mathrm{UTB}+\{c>n:n\in \omega \}.$%
\medskip

\noindent \textbf{(c)} $\mathrm{CT}=\mathrm{PA}(T)+\theta ,$ where $\theta $
is a single sentence that stipulates that $T$ satisfies Tarski's inductive
clauses for a truth predicate for arithmetical sentences.\medskip

\noindent \textbf{(d)} $\mathrm{CT}(c)=\mathrm{CT}+\{c>n:n\in \omega \}.$%
\medskip

The following proposition is well-known and easy to prove; the nontrivial
direction of part (a) is the right-to-left part, which employs a routine
overspill argument; part (b) follows easily from part (a) and the
definitions involved. The proofs of (c) and (d) are routine but somewhat
laborious.\medskip

\noindent \textbf{2.2.3.}~\textbf{Proposition.}~\textit{The following holds
for every model} $\mathcal{M}$ \textit{of} $\mathrm{PA}$ \textit{of any
cardinality}.$\medskip $

\noindent \textbf{(a)} $\mathcal{M}$ \textit{has an inductive partial
satisfaction class iff} $\mathcal{M}$ \textit{has an expansion to a model of
}$\mathrm{UTB}$.\textbf{\medskip }

\noindent \textbf{(b)} $\mathcal{M}$\textit{\ is nonstandard and has an
inductive partial satisfaction class} \textit{iff} $\mathcal{M}$ \textit{has
an expansion to a model of }$\mathrm{UTB}(c)$\textit{.}\textbf{\medskip }

\noindent \textbf{(c)} $\mathcal{M}$ \textit{has an inductive full
satisfaction class iff} $\mathcal{M}$ \textit{has an expansion to a model of
}$\mathrm{CT}$.\textbf{\medskip }

\noindent \textbf{(d)} $\mathcal{M}$\textit{\ is nonstandard and has an
inductive full satisfaction class} \textit{iff} $\mathcal{M}$ \textit{has an
expansion to a model of }$\mathrm{CT}(c)$\textit{.}\textbf{\medskip }

\noindent The concepts of recursive saturation and satisfaction classes are
intimately tied, as witnessed by the following classical result of Barwise
and Schlipf whose proof invokes the resplendence property of countable
recursively saturated models (for a proof, see Corollary 15.12 of \cite%
{Richard-book}). Note that the right-to-left implication in the
Barwise-Schlipf theorem holds for uncountable models $\mathcal{M}$ as well
(and is proved by a simple overspill argument). However, by Tarski's
undefinability of truth theorem the left-to-right direction fails for
`Kaufmann models' (i.e., recursively saturated rather classless
models).\medskip

\noindent \textbf{2.2.4.}~\textbf{Theorem.}~(Barwise-Schlipf) \textit{A
countable model }$\mathcal{M}$ \textit{of }$\mathrm{PA}$\textit{\ is
recursively saturated iff }$\mathcal{M}$\textit{\ has an inductive partial
satisfaction class.}$\medskip $

\begin{center}
\textbf{2.3.}~\textbf{Indiscernibles\medskip }
\end{center}

\noindent \textbf{2.3.1.} \textbf{Definition.}~Given a linear order $(X,<)$,
and nonzero $n\in \omega $, we use $[X]^{n}$ to denote the set of all
\textit{increasing} sequences $x_{1}<\cdot \cdot \cdot <x_{n}$ from $X$.$%
\medskip $

\noindent \textbf{2.3.2.} \textbf{Definition.}~Given a structure $\mathcal{M}
$, some linear order $(I,<)$ where $I\subseteq M$, we say that $(I,<)$ is%
\textit{\ a set of order indiscernibles in} $\mathcal{M}$ if for any $%
\mathcal{L}(\mathcal{M})$-formula $\varphi (x_{1},\cdot \cdot \cdot ,x_{n})$%
, and any two $n$-tuples $\overline{i}$ and $\overline{j}$ from $[I]^{n}$,
we have:

\begin{center}
$\mathcal{M}\models \varphi (i_{1},\cdot \cdot \cdot ,i_{n})\leftrightarrow
\varphi (j_{1},\cdot \cdot \cdot ,j_{n}).$
\end{center}

\noindent \textbf{2.3.3.} \textbf{Definition.}~Suppose $\mathcal{M}$ has
parameter-free definable Skolem functions, and $(I,<_{I})$ is a set of order
indiscernibles in $\mathcal{M}$, and $I_{0}\subseteq I$. We use the notation%
\textit{\ }$\mathcal{M}_{I_{0}}$ \textit{to denote the elementary submodel of%
} $\mathcal{M}$ \textit{generated by} $I_{0}$ (via the parameter-free
definable functions of $\mathcal{M}$).

\begin{itemize}
\item Note that the universe $M_{I_{0}}$ of $\mathcal{M}_{I_{0}}$ consists
of the elements of $M$ that are pointwise definable in $(M,i)_{i\in
I_{0}}.\medskip $
\end{itemize}

\begin{center}
\textbf{2.4.}~\textbf{Interpretability}\medskip
\end{center}

\noindent \textbf{2.4.1.~Definition.~}Suppose $U$ and $V$ are first order
theories, and for the sake of notational simplicity, let us assume that $U$
and $V$ are theories that \textit{support a definable pairing function}. We
use $\mathcal{L}_{U}$ and $\mathcal{L}_{V}$ to respectively designate the
languages of $U$ and $V$.\medskip

\noindent \textbf{(a) }An interpretation $\mathcal{I}$ of $U$ in $V$,
written:

\begin{center}
$U\trianglelefteq _{\mathcal{I}}V$,
\end{center}

\noindent is given by a translation $\tau $ of each $\mathcal{L}_{U}$%
-formula $\varphi $ into an $\mathcal{L}_{V}$-formula $\varphi ^{\tau }$
with the requirement that $V\vdash \varphi ^{\tau }$ for each $\varphi \in U$%
, where $\tau $ is determined by an $\mathcal{L}_{V}$-formula $\delta (x)$
(referred to as a \textit{domain formula}), and a mapping $P\mapsto _{\tau
}A_{P}$ that translates each $n$-ary $\mathcal{L}_{U}$-predicate $P$ into
some $n$-ary $\mathcal{L}_{V}$-formula $A_{P}$. The translation is then
lifted to the full first order language in the obvious way by making it
commute with propositional connectives, and subject to the following clauses:

\begin{center}
$\left( \forall x\varphi \right) ^{\tau }=\forall x(\delta (x)\rightarrow
\varphi ^{\tau })$ and $\left( \exists x\varphi \right) ^{\tau }=\exists
x(\delta (x)\wedge \varphi ^{\tau }).$
\end{center}

\noindent Note that each interpretation $U\trianglelefteq _{\mathcal{I}}V$%
\/gives rise to an \textit{inner model construction} that\textit{\ uniformly
}builds a model $\mathcal{I(M)}\models U$\/ for any $\mathcal{M}\models V$.
\medskip

\noindent \textbf{(b) }$U$ is \textit{interpretable} in $V$ (equivalently: $%
V $\textit{\ interprets} $U$), written $U\trianglelefteq V$, iff $%
U\trianglelefteq _{\mathcal{I}}V$ for some interpretation $\mathcal{I}.$
\medskip

\noindent \textbf{(c) }Given arithmetical theories $U$ and $V$, $U$ is $%
\omega $-\textit{interpretable in }$V$\textit{\ if }the interpretation of
`numbers' and the arithmetical operations of the interpreted theory $U$ are
the same as those of the interpreting theory $V$.\medskip

\noindent \textbf{(d) }$U$ is \textit{locally} \textit{interpretable} in $V$%
, written $U\trianglelefteq _{\mathrm{loc}}V$ if $U_{0}\trianglelefteq V$
for every finite subtheory $U_{0\text{ }}$of $U$.\medskip

\noindent \textbf{(e) }$U$\textbf{\ }and $V$ are \textit{mutually
interpretable} when $U\trianglelefteq V$ and $V\trianglelefteq U.$\medskip

\noindent \textbf{(f) }$U$ is a \textit{retract} of $V$ iff there are
interpretations $\mathcal{I}$ and $\mathcal{J}$ with $U\trianglelefteq _{%
\mathcal{I}}V$ and $V\trianglelefteq _{\mathcal{J}}U$, and a binary $U$%
-formula $F$ such that $F$\ is, $U$-verifiably, an isomorphism between
\textrm{id}$_{U}$ (the identity interpretation on $U$) and $\mathcal{J}\circ
\mathcal{I}$\textsf{.} In model-theoretic terms, this translates to the
requirement that the following holds for every $\mathcal{M}\models U$:

\begin{center}
$F^{\mathcal{M}}:\mathcal{M}\overset{\cong }{\longrightarrow }\mathcal{M}%
^{\ast }:=\mathcal{I(}\mathcal{J(M))}.$
\end{center}

\noindent \textbf{(g) }$U$ and $V$ are \textit{bi-interpretable} iff there
are interpretations $\mathcal{I}$ and $\mathcal{J}$ as above that witness
that $U$ is a retract of $V$, and additionally, there is a $V$-formula $G,$
such that $G$\ is, $V$-verifiably, an isomorphism between the ambient model
of $V$ and the model of $V$ given by $\mathcal{I}\circ \mathcal{J}.$ In
particular, if $U$ and $V$ are bi-interpretable, then given $\mathcal{M}%
\models U$ and $\mathcal{N}\models V$, we have

\begin{center}
$F^{\mathcal{M}}:\mathcal{M}\overset{\cong }{\longrightarrow }\mathcal{M}%
^{\ast }:=\mathcal{I(}\mathcal{J(M))}$ and $G^{\mathcal{N}}:\mathcal{N}%
\overset{\cong }{\longrightarrow }\mathcal{N}^{\ast }:=\mathcal{J(I}\mathcal{%
(N))}.$
\end{center}

\noindent \textbf{(h) }The above notions can also be localized at a pair of
models; in particular suppose $\mathcal{N}$ is an $\mathcal{L}_{U}$%
-structure and $\mathcal{M}$ is an $\mathcal{L}_{V}$-structure. For example,
we say that $\mathcal{N}$ is \textit{interpretable} in $\mathcal{M}$,
written $\mathcal{N}\trianglelefteq \mathcal{M}$ (equivalently: $\mathcal{M}%
\trianglerighteq \mathcal{N}$) iff the universe of discourse of $\mathcal{N}$%
, as well as all the $\mathcal{N}$-interpretations of $\mathcal{L}_{U}$%
-predicates are $\mathcal{M}$-definable. Similarly, we say that $\mathcal{M}$
and $\mathcal{N}$ are \textit{bi-interpretable }if there are parametric
interpretations\textit{\ }$\mathcal{I}$\textit{\ }and\textit{\ }$\mathcal{J}$%
\textit{, }together with an $\mathcal{M}$-definable $F$ and an $\mathcal{N}$%
-definable map $G$ such that:

\begin{center}
$F^{\mathcal{M}}:\mathcal{M}\overset{\cong }{\longrightarrow }\mathcal{M}%
^{\ast }:=\mathcal{I(}\mathcal{J(M))}$ and $G^{\mathcal{N}}:\mathcal{N}%
\overset{\cong }{\longrightarrow }\mathcal{N}^{\ast }:=\mathcal{J(I}\mathcal{%
(M))}.$\medskip
\end{center}

\begin{itemize}
\item Recall that a theory $U$ (with sufficient coding apparatus) is \textit{%
reflexive} if the formal consistency of each finite fragment of $U$ is
provable in $U$.
\end{itemize}

\noindent \textbf{2.4.2.~Theorem.~}(Mostowski's Reflection Theorem) \textit{%
For all }$\mathcal{L}\supseteq \mathcal{L}_{\mathrm{A}}$, \textit{every
extension} (\textit{in the same language}) \textit{of}\textbf{\ }\textrm{PA(}%
$\mathcal{L}$\textrm{)} \textit{is reflexive.}\medskip

\noindent \textbf{2.4.3.~Theorem.~}(Orey's Compactness Theorem)\textbf{\ }%
\textit{If} $U$ \textit{is reflexive, and} $V\trianglelefteq _{\mathrm{loc}%
}U $\textit{\ for some recursively enumerable theory }$V$\textit{, then }$%
V\trianglelefteq U$.\bigskip

\begin{center}
\textbf{3.~THE BASICS OF PAI}\bigskip
\end{center}

\noindent \textbf{3.1.~Definition.~}$\mathrm{PAI}$ is the theory formulated
in $\mathcal{L}_{\mathrm{A}}(I)$ whose axioms are (1) through (3) below.
Note that we write $x\in I$ instead of $I(x)$ for better readability.\medskip

\noindent $(1)$ $\mathrm{PA(}I).$\medskip

\noindent $(2)$ The sentence $\mathrm{Ubd}(I)$ that expresses:
\textquotedblleft $I$\ is unbounded\textquotedblright .\medskip

\noindent $(3)$ The scheme $\mathrm{Indis}_{\mathcal{L}_{\mathrm{A}}}(I)=\{%
\mathrm{Indis}_{\varphi }(I):$ $\varphi $ is an $\mathcal{L}_{\mathrm{A}}$%
-formula$\}$ stipulating that $I$ forms a class of order indiscernibles for
the ambient model of arithmetic$\mathrm{.}$ More explicitly, for each $n$%
-ary formula $\varphi (v_{1},\cdot \cdot \cdot ,v_{n})$ in the language $%
\mathcal{L}_{\mathrm{A}},$ $\mathrm{Indis}_{\varphi }(I)$ is the following
sentence:\medskip

\begin{center}
$\forall x_{1}\in I\cdot \cdot \cdot \forall x_{n}\in I$ $\forall y_{1}\in
I\cdot \cdot \cdot \forall y_{n}\in I$ $[(x_{1}<\cdot \cdot \cdot
<x_{n})\wedge (y_{1}<\cdot \cdot \cdot <y_{n})\rightarrow (\varphi
(x_{1},\cdot \cdot \cdot ,x_{n})\leftrightarrow \varphi (y_{1},\cdot \cdot
\cdot ,y_{n}))].$
\end{center}

\noindent $\mathrm{PAI}^{\circ }$ is the weakening of $\mathrm{PAI}$ in
which the scheme $\mathrm{Indis}_{\mathcal{L}_{\mathrm{A}}}(I)$ is weakened
to the scheme $\mathrm{Indis}_{\mathcal{L}_{\mathrm{A}}}^{\circ }(I)=\{%
\mathrm{Indis}_{\varphi }^{\circ }(I):$ $\varphi $ is an $\mathcal{L}_{%
\mathrm{A}}$-formula$\},$ where $\mathrm{Indis}_{\varphi }^{\circ }(I)$ is
the following sentence:

\begin{center}
$\forall x_{1}\in I\cdot \cdot \cdot \forall x_{n}\in I$ $\forall y_{1}\in
I\cdot \cdot \cdot \forall y_{n}\in I$ $[(x_{1}<\cdot \cdot \cdot
<x_{n})\wedge (y_{1}<\cdot \cdot \cdot <y_{n})\wedge \left( \ulcorner
\varphi \urcorner <x_{1}\wedge \ulcorner \varphi \urcorner <y_{1}\right) $

$\rightarrow (\varphi (x_{1},\cdot \cdot \cdot ,x_{n})\leftrightarrow
\varphi (y_{1},\cdot \cdot \cdot ,y_{n}))].$
\end{center}

\noindent \textbf{3.2.~Proposition.~}\textit{Let }$\mathbb{N}$\textit{\ be
the standard model of }$\mathit{\mathrm{PA.}}$\medskip

\noindent \textbf{(a)} $\mathbb{N}$ \textit{does not have an expansion to a
model of }$\mathrm{PAI}$ (\textit{equivalently: Every model of} $\mathrm{PAI}
$\textit{\ is nonstandard}). \medskip

\noindent \textbf{(b)} $\mathbb{N}$ \textit{has} \textit{an expansion to} $%
\mathrm{PAI}^{\circ }.$\medskip

\noindent \textbf{(c)} \textit{If} $(\mathcal{M},I)$ \textit{is a
nonstandard model of }$\mathrm{PAI}^{\circ }$,\textit{\ and }$c$ \textit{is
any nonstandard element of} $\mathcal{M}$, \textit{then} $(\mathcal{M}%
,I^{>c})\models \mathrm{PAI}$, \textit{where} $I^{>c}=\{i\in I:i>c\}$%
.\medskip

\noindent \textbf{Proof}.\textbf{~}(a)\textbf{\ }is an immediate consequence
of the fact that the standard model of $\mathrm{PA}$ is pointwise definable,
and therefore it does not even have a distinct pair of indiscernibles.
\medskip

To see that (b) holds, fix some enumeration $\left\langle \varphi _{n}:n\in
\omega \right\rangle $ of all arithmetical formulae, and use Ramsey's
theorem to construct a sequence $\left\langle H_{n}:n\in \omega
\right\rangle $ of subsets of $\omega $ such that for each $n\in \omega $
the following three conditions hold:\medskip

\noindent (1) $H_{n}$ is infinite. \medskip

\noindent (2) $H_{n}\supseteq H_{n+1}.$ \medskip

\noindent (3) $H_{n}$ is $\varphi _{n}$-indiscernible (i.e., $\mathrm{Indis}%
_{\varphi _{n}}(H_{n})$ holds).\medskip

\noindent Then recursively define $\left\langle i_{n}:n\in \omega
\right\rangle $ by: $i_{0}=\max \{\min \{H_{0}\},\ulcorner \varphi
_{0}\urcorner \},$ and $i_{n+1}$ is the least $i\in H_{n+1}$ that is greater
than both $i_{n}$ and $\ulcorner \varphi _{n}\urcorner .$ It is easy to see
that $(\mathcal{M},I^{\prime })\models \mathrm{PAI}^{\circ }$, where $%
I=\left\langle i_{n}:n\in \omega \right\rangle $.\footnote{%
The construction of $I$ uses similar ideas as the construction of
(Gaifman's) `minimal types', as in \cite[Theorem 3.1.2]{Roman-Jim's book}}
\medskip

Since (c) readily follows from the definitions involved, the proof is
complete.\hfill $\square $\medskip

\begin{itemize}
\item In light of part (c) of Theorem 3.2, most results in this paper about $%
\mathrm{PAI}$ has a minor variant in which $\mathrm{PAI}$\ is replaced by $%
\mathrm{PAI}^{\circ }.$
\end{itemize}

\noindent \textbf{3.3.~Theorem.~}\textit{Each finite subtheory of }$\mathrm{%
PAI}$ \textit{has an }$\omega $-\textit{interpretation} \textit{in} $\mathrm{%
PA}$. \textit{Consequently}:\medskip

\noindent \textbf{(a) }$\mathrm{PAI}$ \textit{is a conservative extension of}
$\mathrm{PA}$.\medskip

\noindent \textbf{(b) }\textrm{PAI}\ \textit{is interpretable in} \textrm{PA}%
, \textit{hence} \textrm{PA}\textit{\ and }\textrm{PAI} \textit{are mutually
interpretable.}\medskip

\noindent \textbf{(c)}\textit{\ }\textrm{PAI }\textit{is interpretable in}
\textrm{ACA}$_{0}$ (\textit{but not vice-versa}).\footnote{\textrm{ACA}$_{0}$
is the well-known finitely axiomatizable subsystem of second arithmetic that
is conservative over $\mathrm{PA}$.} \medskip

\noindent \textbf{Proof.~}The $\omega $-interpretability of any finite
subtheory of\textit{\ }$\mathrm{PAI}$\ in \textrm{PA}\textit{\ }is an
immediate consequence of the well-known schematic provability of Ramsey's
theorem $\omega \rightarrow (\omega )_{2}^{n}$ in $\mathrm{PA}$ for all
metatheoretic $n\geq 2$ \cite[Theorem 1.5, Chapter II]{Petr-Pavel}. This
makes it evident that (a) holds, and together with Orey's Compactness
Theorem 2.4.3, yields (b). Finally, (c) follows from (b) since $\mathrm{PA}$
is trivially interpretable in \textrm{ACA}$_{0}.$ The parenthetical clause
of (c) is an immediate consequence of (b) and the classical fact that $%
\mathrm{PA}$ is not interpretable in \textrm{ACA}$_{0}$ (the ingredients
whose proof are Mostowski's reflection theorem for \textrm{PA}, finite
axiomatizability of \textrm{ACA}$_{0},$ and G\"{o}del's second
incompleteness theorem). $\square $\medskip

\noindent \textbf{3.4.~Remark.~}Standard techniques can be used to show that
the proof of Theorem 3.3(b)\ yields a feasible reduction of\textit{\ }$%
\mathrm{PAI}$\ in\textit{\ }$\mathrm{PA}$\textit{. }In other words, there is
a polynomial-time function $f$\ such that, given the (binary code)\ of a
proof $\pi $ of an arithmetical sentence $\varphi $ in $\mathrm{PAI}$, $%
f(\pi )$ is the (the binary code of)\ a proof $f(\pi )$ of $\varphi $ in%
\textit{\ }$\mathrm{PA}$. In particular, $\mathrm{PAI}$\ has at most
polynomial speed-up over PA.

\begin{itemize}
\item In what follows $\mathrm{Form}_{k}$ is the set of $\mathcal{L}_{%
\mathrm{A}}$-formulae with precisely $k$ free variables.
\end{itemize}

\noindent \textbf{3.5.~Theorem.~}\textit{The following schemes are provable
in}\textbf{\ }$\mathrm{PAI}$:\ \medskip

\noindent \textbf{(a) }\textit{The \textbf{apartness} scheme:}

\begin{center}
\textit{\ }$\mathrm{\{Apart}_{\varphi }:\varphi \in \mathrm{Form}_{n+1}$, $%
n\in \omega \mathrm{\},}$
\end{center}

\noindent \textit{where} $\mathrm{Apart}_{\varphi }$ \textit{is the
following formula}:\medskip

\begin{center}
$\forall i\in I\ \forall j\in I\left[ i<j\rightarrow \forall x_{1},\cdot
\cdot \cdot ,x_{n}<i\ \left( \exists y\ \varphi (\overline{x},y)\rightarrow
\exists y<j\ \varphi (\overline{x},y)\right) \right] .$
\end{center}

\noindent \textbf{(b)} \textit{The \textbf{diagonal indiscernibility}}%
\footnote{%
This stronger notion of indiscernibility appears often in expositions of the
Paris-Harrington independence result; the same notion is dubbed
\textquotedblleft strong indiscernibility\textquotedblright\ in \cite[%
Definition 3.2.8]{Roman-Jim's book}.}\textit{\ scheme}:

\begin{center}
$\{\mathrm{Indis}_{\varphi }^{+}:$ $\varphi \in \mathrm{Form}_{n+1+r},\
n,r\in \omega ,\ r\geq 1\}$,
\end{center}

\noindent where $\mathrm{Indis}_{\varphi }^{+}(I)$ \textit{is the following
formula}:\medskip

\begin{center}
$\forall i\in I\ \forall \overline{j}\in \lbrack I]^{r}\ \forall \overline{k}%
\in \lbrack I]^{r}\ \left[ \left( i<j_{1}\right) \wedge (i<k_{1})\right]
\longrightarrow $\medskip

$\left[ \forall x_{1},\cdot \cdot \cdot ,x_{n}<i\ \left( \varphi (\overline{x%
},i,j_{1},\cdot \cdot \cdot ,j_{r})\leftrightarrow \varphi (\overline{x}%
,i,k_{1},\cdot \cdot \cdot ,k_{r})\right) \right] .$\medskip
\end{center}

\noindent \textbf{Proof.}~Let $(\mathcal{M},I)\models \mathrm{PAI.}$ To
verify that the apartness scheme holds in $(\mathcal{M},I)$, fix some $%
i_{0}\in I$ and some $\varphi (\overline{x},y)\in \mathrm{Form}_{n+1}.$
Then, since the $I$ is unbounded and the collection scheme holds in $(%
\mathcal{M},I),$ and $I$ is unbounded in $\mathcal{M}$, there is some $%
j_{0}\in I$ with $i_{0}<j_{0}$ such that:

\begin{center}
$(\mathcal{M},I)\models \forall \overline{x}\in \lbrack i_{0}]^{n}\ \left(
\exists y\varphi (\overline{x},y)\rightarrow \exists y<j_{0})\ \varphi (%
\overline{x},y)\right) .$
\end{center}

\noindent The above, together with the indiscernibility of $I$ in $\mathcal{M%
}$, makes it evident that $(\mathcal{M},I)\models \mathrm{Apart}_{\varphi }.$
\medskip

To verify that $\mathrm{Indis}_{\varphi }^{+}(I)$ holds in $(\mathcal{M},I)$%
, we will first establish a weaker form of diagonal indiscernibility of $I$
in which all $j_{n}<k_{1}$ (thus all the elements of $\overline{j}$ are less
than all the elements of $\overline{k}$)$.$ Fix some $\varphi \in \mathrm{%
Form}_{n+1+r}$ and $i_{0}\in I.$ Within $\mathcal{M}$ consider the function $%
f:\left[ M\right] ^{r}\rightarrow \mathcal{P}([i_{0}]^{n})$ by:

\begin{center}
$f(\overline{y}):=\{\overline{a}\in \lbrack i_{0}]^{n}:\varphi (\overline{a}%
,i_{0},\overline{y})\}.$
\end{center}

\noindent Since $(\mathcal{M},I)$ satisfies the collection scheme and $I$ is
unbounded in $\mathcal{M}$, this shows there are $y_{1}<\cdot \cdot \cdot
<y_{2r}$ in $I$ such that:

\begin{center}
$f(y_{1},\cdot \cdot \cdot ,y_{r})=f(y_{r+1},\cdot \cdot \cdot ,y_{2r}).$
\end{center}

\noindent Thus $(\mathcal{M},I)$ satisfies:

\begin{center}
$\forall \overline{x}\in \lbrack i_{0}]^{n}\ \left[ \varphi (\overline{x}%
,i_{0},y_{1},\cdot \cdot \cdot ,y_{r})\leftrightarrow \varphi (\overline{x}%
,i_{0},y_{r+1},\cdot \cdot \cdot ,y_{2r})\right] .$
\end{center}

\noindent By the indiscernibility of $I$ in $\mathcal{M}$, the above implies
the following weaker form of $\mathrm{Indis}_{\varphi }^{+}(I)$:\medskip

\begin{center}
$\forall i\in I\ \forall \overline{j}\in \lbrack I]^{r}\ \forall \overline{k}%
\in \lbrack I]^{r}\ \left[ \left( i<j_{1}\right) \wedge (j_{n}<k_{1})\right]
\longrightarrow $\medskip

$\left[ \forall \overline{x}\in (i_{0})^{n}\ \left( \varphi (\overline{x}%
,i,j_{1},\cdot \cdot \cdot ,j_{r})\leftrightarrow \varphi (\overline{x}%
,i,k_{1},\cdot \cdot \cdot ,k_{r})\right) \right] .$\medskip
\end{center}

\noindent We will now show that the above weaker form of $\mathrm{Indis}%
_{\varphi }^{+}(I)$ already implies $\mathrm{Indis}_{\varphi }^{+}(I).$
Given $i\in I$, $\overline{a}\in \lbrack I]^{r}\ $and $\overline{b}\in
\lbrack I]^{r},$ with $i<a_{1}$ and $i<b_{1},$ choose $\overline{y}\in
\lbrack I]^{r}$ with $y_{1}>\max \left\{ a_{n},b_{n}\right\} .$ Then by the
above we have:

\begin{center}
$\mathcal{M}\models \left[ \forall \overline{x}\in (i)^{n}\ \left( \varphi (%
\overline{x},i,a_{1},\cdot \cdot \cdot ,a_{r})\leftrightarrow \varphi (%
\overline{x},i,y_{1},\cdot \cdot \cdot ,y_{r})\right) \right] ,$
\end{center}

\noindent and

\begin{center}
$\mathcal{M}\models \left[ \forall \overline{x}\in (i)^{n}\ \left( \varphi (%
\overline{x},i,b_{1},\cdot \cdot \cdot ,b_{r})\leftrightarrow \varphi (%
\overline{x},i,y_{1},\cdot \cdot \cdot ,y_{r})\right) \right] ,$
\end{center}

\noindent which together imply:

\begin{center}
$\mathcal{M}\models \left[ \forall \overline{x}\in (i)^{n}\ \left( \varphi (%
\overline{x},i,a_{1},\cdot \cdot \cdot ,a_{r})\leftrightarrow \varphi (%
\overline{x},i,b_{1},\cdot \cdot \cdot ,b_{r})\right) \right] .$

\hfill $\square $\medskip
\end{center}

\begin{itemize}
\item Note that the diagonal indiscernibility scheme for $\mathcal{L}_{%
\mathrm{A}}$-formulae ensures that if $(\mathcal{M},I)\models \mathrm{PAI}$
and $i\in I,$ then $I^{\geq i}$ is a set of indiscernibles over the expanded
structure $(\mathcal{M},m)_{m<i}$, where $I^{\geq i}=\{x\in I:x\geq i\}.$
\end{itemize}

\pagebreak

\begin{center}
\textbf{4.~MAIN RESULTS}\bigskip
\end{center}

In this section we prove refinements of Theorems A and B of the abstract (as
in Theorems 4.6 and 4.12).\medskip

\noindent \textbf{4.1.~Theorem.~}\textit{There is a formula }$\sigma (x)$
\textit{in the language} $\mathcal{L}_{\mathrm{A}}(I)$ \textit{such that }$%
S=\sigma ^{\mathcal{M}}$ \textit{is an inductive partial satisfaction class}
\textit{on} $\mathcal{M}$ \textit{for all models} $(\mathcal{M},I)\models
\mathrm{PAI.}$ \medskip

\noindent \textbf{Proof.}~We first define a recursive function that
transforms each formula $\varphi (\overline{x})\in \mathrm{Form}_{n}$ into a
$\Delta _{0}$-formula $\varphi ^{\ast }(\overline{x},z_{1},\cdot \cdot \cdot
,z_{k})$, where $\{z_{i}:1\leq i\in \omega \}$ is a fresh supply of
variables added to the syntax of first order logic (the definition of $%
\varphi ^{\ast }$ below will make it clear that $k$ is the $\exists $-depth
of $\varphi )$. In what follows $x$ and $y$ range over the set of variables
before the addition of the fresh stock of $z_{i}$s$.$ We assume that the
only logical constants used in $\varphi $ are $\left\{ \lnot ,\vee ,\exists
\right\} $ and none of the fresh variables $z_{i}$ occurs in $\varphi
.\medskip $

\noindent $(1)$ If $\varphi $ is an atomic $\mathcal{L}_{\mathrm{A}}$%
-formula, then $\varphi ^{\ast }=\varphi $.\medskip

\noindent $(2)$ $\left( \lnot \varphi \right) ^{\ast }=\lnot \varphi ^{\ast
}.$\medskip

\noindent $(3)$ $(\varphi _{1}\vee \varphi _{2})^{\ast }=\varphi _{1}^{\ast
}\vee \varphi _{2}^{\ast }.$\medskip

\noindent $(4)$ $\left( \exists y\ \varphi \right) ^{\ast }=\exists y<z_{1}$
$\widetilde{\varphi ^{\ast }},$ where $\varphi ^{\ast }=\varphi ^{\ast }(%
\overline{x},y,z_{1},\cdot \cdot \cdot ,z_{k}),$ and $\widetilde{\varphi
^{\ast }}$ is the result of replacing $z_{i}$ with $z_{i+1}$ in $\varphi
^{\ast }$ for each $1\leq i\leq k.$ \medskip

\noindent \textbf{Claim }$\mathbf{(}\nabla )$\textbf{.}~\textit{Suppose }$%
\varphi =\varphi (\overline{x})\in \mathrm{Form}_{n}$, and $\varphi ^{\ast
}=\varphi ^{\ast }(\overline{x},z_{1},\cdot \cdot \cdot ,z_{k}),$ $(\mathcal{%
M},I)\models \mathrm{PAI}$, $\overline{a}\in M^{n}$, \textit{and} $%
(i_{1},\cdot \cdot \cdot ,i_{k})\in \lbrack I]^{k}$ \textit{such that there
is some} $j\in I$ \textit{with} $j<i_{1}$ \textit{and} $a_{s}<j$ \textit{for
each} $1\leq s\leq n.$ \textit{Then }$\mathcal{M}$ \textit{satisfies}:

\begin{center}
$\varphi (\overline{a})\leftrightarrow \varphi ^{\ast }(\overline{a}%
,i_{1},\cdot \cdot \cdot ,i_{k})$.
\end{center}

\noindent \textbf{Proof.}~We use induction of the complexity of $\varphi .$
The only case that needs an explanation is the existential case, the others
go through trivially. Thus, it suffices to verify that if $(i_{1},\cdot
\cdot \cdot ,i_{k+1})\in \lbrack I]^{k}$ and there is some $j\in I$ with $%
j<i_{1}$ and $a_{s}<j$ for each $1\leq s\leq k,$ then:\medskip

\noindent $(\nabla )$\ \ \ $\mathcal{M}\models \left( \exists y\ \varphi (%
\overline{a},y)\leftrightarrow \exists y<i_{1}\ \varphi ^{\ast }(\overline{a}%
,y,i_{2},\cdot \cdot \cdot ,i_{k+1})\right) $,\medskip

\noindent where $\left( \varphi (\overline{x},y)\right) ^{\ast }=\varphi
^{\ast }(\overline{x},y,z_{1},\cdot \cdot \cdot ,z_{k}).$ To establish the
left-to-right direction of $(\nabla )$, suppose $\mathcal{M}\models \exists
y\ \varphi (\overline{a},y).$ By the veracity of the apartness scheme and
the assumption that $a_{s}<j$ for each $1\leq s\leq n$, there is $b<i_{1}$
such that $\mathcal{M}\models \varphi (\overline{a},b)$. Thus since $b$, as
well as $a_{1},\cdot \cdot \cdot ,a_{n}$ are all below $i_{1}$, $i_{1}$ can
serve as the element \textquotedblleft $j$\textquotedblright\ of the
inductive assumption, hence allowing us to conclude that $\mathcal{M}\models
\varphi (\overline{a},b)$ iff $\varphi ^{\ast }(\overline{a},b,i_{2},\cdot
\cdot \cdot ,i_{k+1}),$ therefore $\mathcal{M}\models \exists y<i_{1}$ $%
\varphi ^{\ast }(\overline{a},y,i_{2},\cdot \cdot \cdot ,i_{k+1}),$ as
desired. The right-to-left direction of $(\nabla )$ is trivial. This
concludes the proof of the claim $\mathbf{(}\nabla )$.\medskip

We are now ready to show that there is an $(\mathcal{M},I)$-definable $%
S\subseteq M$ such that $S$ is an inductive satisfaction class over $%
\mathcal{M}$\textit{.} The following procedure takes place in $(\mathcal{M}%
,I),$ in particular, the variables $n$ and $k$ in ($\mathbb{P}$) range in $M$
and need not be standard: \medskip

\noindent ($\mathbb{P}$)\ \ \ Given any $\varphi (\overline{x})\in \mathrm{%
Form}_{n}$ and any $n$-tuple $\overline{a}$, calculate $\left( \varphi (%
\overline{x})\right) ^{\ast }=\varphi ^{\ast }(\overline{x},z_{1},\cdot
\cdot \cdot ,z_{k}),$ and let $j\in I$ be the first element of $I$ such that
$\ulcorner \varphi (\overline{x})\urcorner <j$ and $a_{s}<j$ for each $1\leq
s\leq n$, and then let and $i_{1},\cdot \cdot \cdot ,i_{k}$ to be the first $%
k$ elements of $I$ that are above $j.$ Then define $S$ by:

\begin{center}
$\left\langle \varphi ,\overline{a}\right\rangle \in S$ iff $[\varphi ^{\ast
}(\overline{a},i_{1},\cdot \cdot \cdot ,i_{k})\in \mathrm{Sat}_{\Delta
_{0}}],$
\end{center}

\noindent where $\mathrm{Sat}_{\Delta _{0}}$ is the canonical $\Sigma _{1}$%
-definable satisfaction predicate for $\Delta _{0}$ formulae of arithmetic.
\medskip

\noindent Thus the desired formula $\sigma $ is given by

\begin{center}
$\sigma (\left\langle \varphi ,\overline{a}\right\rangle ):=[\varphi ^{\ast
}(\overline{a},i_{1},\cdot \cdot \cdot ,i_{k})\in \mathrm{Sat}_{\Delta
_{0}}].$
\end{center}

\hfill $\square $\medskip

\noindent \textbf{4.2.~Remark.}~Three remarks are in order concerning the
proof of Theorem 4.1.\medskip

\noindent \textbf{(a)} If $\varphi (\overline{x})$ is a standard formula, $(%
\mathcal{M},I)\models \mathrm{PAI}$, and $j\in I,$ then the condition $%
\ulcorner \varphi (\overline{x})\urcorner <j$ in the procedure ($\mathbb{P)}$
is automatically satisfied since every element of $I$ is nonstandard. The
role of the condition $\ulcorner \varphi (\overline{x})\urcorner <j$ will
become clear in the proof of Lemma 4.10 and Theorem 4.11.\medskip

\noindent \textbf{(b)} If $I^{\prime }$ is a cofinal subset of $I$ such that
$(\mathcal{M},I^{\prime })\models \mathrm{PAI}$ and $S^{\prime }$ is the
partial satisfaction class on $\mathcal{M}$ as defined by $\sigma $ in $(%
\mathcal{M},I^{\prime }),$ then thanks to the diagonal indiscernibility
property of $I$, $S=S^{\prime }.$ This fact comes handy in the proof of
Theorem 5.5.\medskip

\noindent \textbf{(c) }The transformation $\varphi \mapsto \varphi ^{\ast }$
given in the proof of Theorem 4.1 can be reformulated in the following more
intuitive way: Given $\varphi (\overline{x})\in \mathrm{Form}_{n}$, find an
equivalent formula $\varphi ^{\mathrm{pnf}}(\overline{x})$ in the prenex
normal form:

\begin{center}
$\varphi ^{\mathrm{pnf}}(\overline{x})=\forall v_{1}\exists w_{1}\cdot \cdot
\cdot \ \delta (v_{1},w_{1}\cdot \cdot \cdot ,v_{k},w_{k},\overline{x})$,
\end{center}

\noindent and then define $\left( \varphi (\overline{x})\right) ^{\ast }$ to
be:

\begin{center}
$\forall v_{1}<z_{1}\ \exists w_{1}<z_{2}\cdot \cdot \cdot \ \delta
(v_{1},w_{1},\cdot \cdot \cdot ,v_{k},w_{k},\overline{x}).$
\end{center}

\noindent A similar transformation is found in the proof of the
Paris-Harrington Theorem \cite{Paris-Harrington}. $\medskip $

\noindent \textbf{4.3.~Corollary}. \textit{The following hold for every
model }$\mathcal{M}$ \textit{of} $\mathrm{PA}$ \textit{of any cardinality}.
\medskip

\noindent \textbf{(a)} \textit{There is no} $\mathcal{M}$-\textit{definable
subset} $I$ \textit{of} $M$\textit{\ such that} $(\mathcal{M},I)\models
\mathrm{PAI}$ (\textit{therefore no rather classless recursively saturated
model of} $\mathrm{PA}$ \textit{has an expansion to a model of} $\mathrm{PAI)%
}$.\medskip

\noindent \textbf{(b)} \textit{If }$\mathcal{M}$ \textit{has an expansion to
a model of} \textrm{PAI}, \textit{then} $\mathcal{M}$\ \textit{is
recursively saturated; and the converse holds if }$\mathcal{M}$\textit{\ is
countable.}\medskip

\noindent \textbf{(c)} \textit{If }$\mathcal{M}$ \textit{has an expansion }$(%
\mathcal{M},I)\models \mathrm{PAI}$, \textit{then} $M\neq M_{I}$, \textit{%
where} $M_{I}$ \textit{consists of elements of} $M$ \textit{that are
definable in} $(\mathcal{M},i)_{i\in I}.$\medskip

\noindent \textbf{Proof.}~(a) follows by putting Theorem 4.1 together with
Tarski's theorem on undefinability of truth.\footnote{%
An alternative, more direct proof of (a) invokes diagonal indiscernibility.
Suppose to the contrary that $(\mathcal{M},I)\models \mathrm{PAI}$ and $I$
is definable in $\mathcal{M}$ by a formula $\varphi (x,m)$ for some $m\in M.$
Let $i_{1}<i_{2}<i_{3}$ be the first three elements of $I$ above $\mathcal{M}
$. Note that $i_{2}$ and $i_{3}$ is each pointwise definable in $(\mathcal{M}%
,m,I).$ Hence $i_{2}$ and $i_{3}$ are discernible in $(\mathcal{M},m,I),$
and therefore they are also discernible in $(\mathcal{M},m)$ (since $I$ is
definable in $\mathcal{M}$ with parameter $m$). On the other hand, by the
diagonal indiscernibility property of $I$, for any arithmetical formula $%
\theta (x,y)$, $\mathcal{M}$ satisfies $\theta (m,i_{1})\leftrightarrow
\theta (m,i_{2}).$ We have arrived at a contradiction.} (b) follows directly
by putting Theorem 4.1 with the Barwise-Schlipf Theorem 2.2.4. To verify (c)
suppose $M_{I}=M$ for $(\mathcal{M},I)\models \mathrm{PAI}$. Recall that $%
\mathcal{M}$ is nonstandard by Proposition 3.2(a). By Theorem 4.1 there is
an inductive partial satisfaction class $S$ on $\mathcal{M}$ that is
definable in $(\mathcal{M},I)$. Consider the function

\begin{center}
$h:M\rightarrow M,$
\end{center}

\noindent where $h$ is defined in $(\mathcal{M},I)$ by $h(m):=$ the (G\"{o}%
del number of) the least $\mathcal{L}_{\mathrm{A}}$-formula $\varphi (x,%
\overline{y})$ such that, as deemed by $S$, $m$ is defined by $\varphi (x,%
\overline{i})$ for some tuple $\overline{i}$ of parameters from $I,$ i.e., $%
S $ contains the sentences $\varphi (m,\overline{i})$ and $\exists !x\
\varphi (x,\overline{i}).$ Note that the set of \textit{standard} elements
of $\mathcal{M}$ is definable in $(\mathcal{M},I)$ as the set of $i$ such
that $i<j$ for some $j$ in the range of $h$. Thus $(\mathcal{M},I)$ is a
nonstandard model of \textrm{PAI}, in which the standard cut $\omega $ is
definable, which is impossible.\hfill $\square $\medskip

\noindent \textbf{4.4}.~\textbf{Remark.}~As shown by Schmerl \cite{Jim
indiscernibles}, every countable recursively saturated model $\mathcal{M}$
of $\mathrm{PA}$ carries a set of indiscernibles $I$ such that $M_{I}=M.$
Thus, in light of part (c) of Corollary 4.3, such a set of indiscernibles $I$
never has the property that $(\mathcal{M},I)\models \mathrm{PAI.}$\medskip

\noindent \textbf{4.5}.~\textbf{Remark.}~In contrast to part (c) of
Corollary 4.3, the proof technique of the Kossak-Schmerl construction of
prime inductive partial satisfaction classes (as in Theorem 10.5.2 of \cite%
{Roman-Jim's book}) can be readily adapted to show that every countable
recursively saturated model $\mathcal{M}$ of $\mathrm{PA}$ has a pointwise
definable expansion $(\mathcal{M},I)\models \mathrm{PAI.}$\medskip

\noindent \textbf{4.6}.~\textbf{Theorem.}~\textit{The following are
equivalent for a model }$\mathcal{M}$ \textit{of }$\mathrm{PA}$\textit{\ of
any cardinality}:$\medskip $

\noindent $(i)$ $\mathcal{M}$ \textit{has an expansion to a model of }$%
\mathrm{UTB}(c)$\textit{.}\medskip

\noindent $(ii)$ $\mathcal{M}$ \textit{has an expansion to a model of }$%
\mathrm{PAI.}$\ \medskip

\noindent \textit{Consequently}, $\mathcal{M}$ \textit{has an expansion to a
model of }$\mathrm{UTB}$ \textit{iff} $\mathcal{M}$ \textit{has an expansion
to a model of }$\mathrm{PAI}^{\circ }\mathrm{.}$\medskip

\noindent \textbf{Proof.} Since $(ii)\Rightarrow (i)$ is justified by
Theorem 4.1, it suffices to show that $(i)\Rightarrow (ii).$\footnote{%
As pointed out by Roman Kossak, $(i)\Rightarrow (ii)$ was first noted in
\cite[Proposition 4.5]{Roman (omega-property)}.} By Proposition 2.2.3(b)
there is an inductive partial satisfaction class $S$ on $\mathcal{M}.$
Consider the $\mathcal{L}_{\mathrm{A}}(S)$-formula $\psi (x)$ that expresses:

\begin{center}
\textquotedblleft there is a definable (in the sense of $S$) unbounded set
of indiscernible for $\mathcal{L}_{\mathrm{A}}$-formulae of G\"{o}del-number
at most $x$\textquotedblright .
\end{center}

\noindent More specifically, $\psi (x)$ is the formula $\exists \theta \in
\mathrm{Form}_{1}$ ($U(\theta )\wedge H(\theta ,x)),$ where $U(\theta )$ is
the following $\mathcal{L}_{\mathrm{A}}(S)$-sentence:

\begin{center}
$\left[ \forall x\exists y(x<y\wedge \theta (x))\right] \in S$,
\end{center}

\noindent and $H(\theta ,x)$ is the following $\mathcal{L}_{\mathrm{A}}(S)$%
-sentence:

\begin{center}
$\forall \varphi \in \mathrm{Form}(\varphi \leq x\rightarrow \mathrm{Indis}%
_{\varphi }(\theta )\in S),$
\end{center}

\noindent where $\mathrm{Indis}_{\varphi }(\theta )$ is the following $%
\mathcal{L}_{\mathrm{A}}$-sentence:

\begin{center}
$\forall x_{1}\cdot \cdot \cdot \forall x_{2n}$ $[(x_{1}<\cdot \cdot \cdot
<x_{n})\wedge (x_{n+1}<\cdot \cdot \cdot <x_{2n})\wedge
\bigwedge\limits_{1\leq i\leq 2n}\theta (x_{i})]\rightarrow $

$(\varphi (x_{1},\cdot \cdot \cdot ,x_{n})\leftrightarrow \varphi
(x_{n+1},\cdot \cdot \cdot ,x_{2n}))].$
\end{center}

\noindent By the schematic provability of Ramsey's theorem in $\mathrm{PA}$,
$(\mathcal{M},S)\models \psi (n)$ for each $n\in \omega $, so by overspill, $%
(\mathcal{M},S)\models \psi (c)$ holds for some nonstandard $c\in M$. Hence
there is some $\theta _{0}\in \mathrm{Form}_{1}^{\mathcal{M}}$ such that $(%
\mathcal{M},S)\models H(\theta ,c)$, and thus $(\mathcal{M},I)\models
\mathrm{PAI}$, where:

\begin{center}
$I:=\{m\in M:(\mathcal{M},S)\models \theta (m)\in S\}$.
\end{center}

\noindent This concludes the proof of the equivalence of $(i)$ and $(ii)$.
The `consequently' clause readily follows from Proposition 3.2 and the
equivalence of $(i)$ and $(ii)$.\hfill $\square $\medskip

\begin{itemize}
\item Going back to Theorem 4.1, one might wonder if it is possible for $%
\sigma ^{\mathcal{M}}$ to be a \textit{full} satisfaction class on $\mathcal{%
M}$. There are certainly many models $(\mathcal{M},I)$ of $\mathrm{PA}$ for
which $\sigma ^{\mathcal{M}}$ is not a full satisfaction class since the
existence of a full inductive satisfaction class on a model $\mathcal{M}$
implies that $\mathrm{Con}(\mathrm{PA})$ holds in $\mathcal{M}$ (and much
more, see the remarks following Theorem 4.9). The results of the rest of
this section are informed by this question.\medskip
\end{itemize}

\noindent \textbf{4.7.~Definition.~}$\alpha $ is the $\mathcal{L}_{\mathrm{A}%
}(I)$-sentence the sentence $\alpha $ expressing \textquotedblleft $\sigma $
defines a full satisfaction class\textquotedblright , where $\sigma (x)$ is
the formula given in the proof of Theorem 4.1.\medskip

\noindent \textbf{4.8.~Definition.~}Given a recursively axiomatized theory $%
T $ extending $\mathrm{I}\Delta _{0}+\mathrm{Exp}$, the \textit{uniform
reflection scheme over} $T$, denoted $\mathrm{RFN}(T)$, is defined via:

\begin{center}
$\mathrm{RFN}(T):=\{\forall x($\textrm{Prov}$_{T}(\ulcorner \varphi (\overset%
{\cdot }{x})\urcorner )\rightarrow \varphi (x)):\varphi (x)\in \mathrm{Form}%
_{\mathrm{1}}\}.$
\end{center}

\noindent The sequence of schemes $\mathrm{RFN}^{\alpha }(T)$, where $\alpha
$ is recursive ordinal $\alpha $, is defined as follows: \medskip

\noindent $\mathrm{RFN}^{0}(T)=T;$\medskip

\noindent $\mathrm{RFN}^{\alpha +1}(T)=\mathrm{RFN}(\mathrm{RFN}^{\alpha
+1}(T));$\medskip

\noindent $\mathrm{RFN}^{\gamma }(T)=\bigcup\limits_{\alpha <\gamma }\mathrm{%
RFN}^{\alpha }(T).$\medskip

\noindent \textbf{4.9.~Theorem.~}(Folklore) \textit{The arithmetical
consequences of} \textrm{CT} \textit{are axiomatized by }$\mathrm{PA}+%
\mathrm{RFN}^{\varepsilon _{0}}(\mathrm{PA}).$\medskip

\noindent \textbf{Proof.~}It is well-known \cite[Section 8.6]{Halbach book}
that the arithmetical consequences of $\mathrm{PA}(S)+\mathrm{FS}(S)$
coincide with the arithmetical consequences of $\mathrm{ACA}$ (the extension
of $\mathrm{ACA}_{0}$ by the full induction scheme). It has long been known
that the arithmetical consequences of $\mathrm{ACA}$ can be axiomatized by $%
\mathrm{PA}+\mathrm{RFN}^{\varepsilon _{0}}(\mathrm{PA}),$ a result which
has been recently given a new proof in the work of Beklemishev and Pakhomov
\cite[Sec.~8.3]{Lev+Fedya}.\footnote{%
It is also known that $\mathrm{PA}+\mathrm{RFN}^{\varepsilon _{0}}(\mathrm{PA%
})$ can be axiomatized by $\mathrm{PA}+\mathrm{TI}(\varepsilon _{\varepsilon
_{0}})$, where $\mathrm{TI}(\varepsilon _{\varepsilon _{0}})$ is the scheme
of transfinite induction for ordinals less than $\varepsilon _{\varepsilon
_{0}}$ ($\varepsilon _{\alpha }$ is the $\alpha $-th $\varepsilon $-number,
i.e., the $\alpha $-th fixed point of the map $\gamma \mapsto \omega
^{\gamma }).$}\hfill $\square $\medskip \medskip

The following lemma, which will come handy at the end of the proof of
Theorem 4.11, shows that if $(\mathcal{M},I)\models \mathrm{PAI}+\alpha $,
and $\varphi \in \mathrm{Form}^{\mathcal{M}}$ (note that $\varphi $ is
allowed to be nonstandard), then as viewed by $S$, a tail of $I$ satisfies
diagonal $\varphi $-indiscernibility.\medskip

\noindent \textbf{4.10.~Lemma.~}\textit{Suppose }$(\mathcal{M},S,I)\models
\mathrm{PAI}+\mathrm{PA(}S,I)+\mathrm{FS}(S).$ \textit{If} $\varphi \in
\mathrm{Form}_{n+r+1}^{\mathcal{M}}$ (\textit{where} $n,r\in M),$ \textit{%
then}:

\begin{center}
$(\mathcal{M},S,I)\models \forall i\in I\ \left( \varphi <i\ \longrightarrow
\theta (S,i,\varphi )\right) ,$
\end{center}

\noindent \textit{where} $\theta (i,\varphi )$ \textit{is the following} $%
\mathcal{L}_{\mathrm{A}}(S,I)$-\textit{formula}:

\begin{center}
$\forall \overline{j}\in \lbrack I]^{r}\ \forall \overline{k}\in \lbrack
I]^{r}\ \left[ \left( i<j_{1}\right) \wedge (i<k_{1})\right] \longrightarrow
$\medskip

$\left[ \forall x_{1},\cdot \cdot \cdot ,x_{n}<i\ \left( \varphi (\overline{x%
},i,j_{1},\cdot \cdot \cdot ,j_{r})\in S\leftrightarrow \varphi (\overline{x}%
,i,k_{1},\cdot \cdot \cdot ,k_{r})\in S\right) \right] .$
\end{center}

\noindent \textbf{Proof.~}The strategy of establishing the diagonal
indiscernibility of $I$ in the proof of Theorem 3.5(b) can be readily
carried out in this context, thanks to the fact that $(\mathcal{M}%
,S,I)\models \mathrm{PA(}S,I)$.\hfill $\square $\medskip

\noindent \textbf{4.11.~Theorem.~}\textit{There is a formula }$\iota (x)$
\textit{in the language} $\mathcal{L}_{\mathrm{A}}(T,c)$ \textit{such that
for all models }$(\mathcal{M},T,c)$\textit{\ of }$\mathrm{CT}(c)$\textit{, }$%
(\mathcal{M},I)\models \mathrm{PAI}+\alpha $\textrm{\ }\textit{for} $I=\iota
^{(\mathcal{M},T,c)}$.\medskip

\noindent \textbf{Proof.~}We will describe the formula\textbf{\ }$\iota (x)$
by working in an arbitrary model $(\mathcal{M},T,c)\models \mathrm{CT}(c).$
Since $\mathrm{PA}(S)+\mathrm{FS}(S)$ and $\mathrm{CT}$ are well-known to be
bi-interpretable, we do most of our work with the model $(\mathcal{M}%
,S)\models \mathrm{PA}(S)+\mathrm{FS}(S)$, and at the end will take
advantage of the nonstandard element $c$. The basic idea is that $\mathrm{PA}
$ can verify that the formalized (infinite)\ Ramsey theorem is provable in $%
\mathrm{PA}$, so using the inductive full satisfaction class $S$ we can
follow the strategy of the proof of part (b) of Theorem 3.2 to define a set $%
I$ in $(\mathcal{M},S)$ such that not only $(\mathcal{M},I)\models \mathrm{%
PAI}+\alpha \mathrm{.}$ More specifically, Ramsey's theorem's can be
fine-tuned by asserting that if an arithmetically definable coloring $f$ of $%
m$-tuples is of complexity $\Sigma _{n}$, then there is an arithmetical
infinite monochromatic subset for $f$ of complexity $\Sigma _{n+m+1}$ \cite%
{Petr-Pavel}. Therefore

\begin{center}
$\mathrm{PA}\vdash \forall r\geq 2\ \forall \varphi \in \mathrm{Form}_{r}\
\exists \theta \in \mathrm{Form}_{1}\ \mathrm{Prov}_{\mathrm{PA}}(\mathrm{%
Indisc}_{\varphi }(\theta )),$
\end{center}

\noindent where $\mathrm{Indisc}_{\varphi }(\theta )$ is as in the proof of
Theorem 4.6. On the other hand, it is well-known that $\mathrm{PA}(S)+%
\mathrm{FS}(S)$ proves the global reflection principle\footnote{%
Indeed, as shown in \L e\l yk's dissertation \cite{Mateusz thesis}, the
global reflection principle can be proved in the fragment \textrm{CT}$_{0}$
of \textrm{CT}.}

\begin{center}
$\forall \varphi ($\textrm{Prov}$_{\mathrm{PA}}(\varphi )\rightarrow
S(\varphi )).$
\end{center}

\noindent Hence \medskip

\noindent $(\ast )\ \ \ (\mathcal{M},S)\models \forall r\geq 2\ \forall
\varphi \in \mathrm{Form}_{r}\ \exists \theta \in \mathrm{Form}_{1}\ \mathrm{%
Indisc}_{\varphi }(\theta )\in S.$\medskip

\noindent Reasoning in $(\mathcal{M},S)$, fix some enumeration $\left\langle
\varphi _{m}:m\in M\right\rangle $ of all arithmetical formulae, and use $%
(\ast )$ to construct a recursive sequence $\left\langle \theta _{m}:m\in
M\right\rangle $ of elements of $\mathrm{Form}_{1}^{\mathcal{M}}$ such that
the following three conditions hold:\medskip

\noindent (1) $(\mathcal{M},S)\models \forall m\forall x\ \exists y>x\
\theta _{m}(\overset{\cdot }{y})\in S.$\medskip

\noindent (2) $(\mathcal{M},S)\models \forall m\forall x(\theta _{m+1}(%
\overset{\cdot }{x})\in S\rightarrow \theta _{m}(\overset{\cdot }{x})\in S).$
\medskip

\noindent (3) $(\mathcal{M},S)\models \forall m\ \mathrm{Indisc}_{\varphi
_{m}}(\theta _{m})\in S$.\medskip

\noindent Let $H_{m}=\{x\in M:(\mathcal{M},S)\models \theta _{m}(\overset{%
\cdot }{x})\in S\}$, and within $(\mathcal{M},S)$ recursively define $%
\left\langle i_{m}:m\in M\right\rangle $ by: $i_{0}=\max \{\min
\{H_{0}\},\ulcorner \varphi _{0}\urcorner \},$ and $i_{m+1}$ is the least $%
i\in H_{m+1}$ that is greater than both $i_{m}$ and $\ulcorner \varphi
_{m}\urcorner .$ It is easy to see that $(\mathcal{M},I)\models \mathrm{PAI}%
^{\circ }$, where $I=\left\langle i_{m}:m\in M\right\rangle ,$ and therefore
as noted in Proposition 3.2(c) $(\mathcal{M},I^{>c})\models \mathrm{PAI}$.
The procedure described for constructing $I$ makes it clear that $I$ is
definable by an $\mathcal{L}_{\mathrm{A}}(T,c)$-formula $\iota (x).$ \medskip

It remains to show that $(\mathcal{M},I)\models \mathrm{PAI}+\alpha .$ Note
that $(\mathcal{M},S,I)\models \mathrm{PA(}S,I),$ this is precisely where
Lemma 4.10 comes to the rescue, since together with the veracity of $\mathrm{%
PAI}+\mathrm{PA(}S,I)+\mathrm{FS}(S)$ in $(\mathcal{M},S,I)$ it allows us
verify the following nonstandard analogue $\mathbf{(}\nabla ^{\ast })$ of $%
\mathbf{(}\nabla )$ from the proof of Theorem 4.1 (in what follows the map $%
\varphi \mapsto \varphi ^{\ast }$ is defined as in the proof of Theorem 4.1
within $\mathcal{M)}$.\medskip

\noindent $\mathbf{(}\nabla ^{\ast })$\textbf{.}~Suppose\textit{\ }$\varphi
=\varphi (\overline{x})\in \mathrm{Form}_{r}^{\mathcal{M}}$ for some $r\in M$
(NB: $r$ need not be standard), $\varphi ^{\ast }=\varphi ^{\ast }(\overline{%
x},z_{1},\cdot \cdot \cdot ,z_{k-1})$, where $k\in M,$ $\overline{a}\in
M^{r} $, \textit{and} $(i_{1},\cdot \cdot \cdot ,i_{k})\in \lbrack I]^{k}$
such that there is some $j\in I$ \textit{with} $j<i_{1}$ and $a_{s}<j$ for
each $1\leq s\leq r.$ Then $(\mathcal{M},S,I)$ satisfies the following:

\begin{center}
$\varphi (\overline{a})\in S\leftrightarrow \lbrack \varphi ^{\ast }(%
\overline{a},i_{1},\cdot \cdot \cdot ,i_{k})\in \mathrm{Sat}_{\Delta _{0}}]$.
\end{center}

\noindent Recall that $\sigma (\left\langle \varphi ,\overline{a}%
\right\rangle ):=[\varphi ^{\ast }(\overline{a},i_{1},\cdot \cdot \cdot
,i_{k})\in \mathrm{Sat}_{\Delta _{0}}].$ $\mathbf{(}\nabla ^{+})$ assures us
that $\sigma ^{\mathcal{M}}$ coincides with $S$, and thus $(\mathcal{M}%
,I)\models \mathrm{PAI}+\alpha $.\hfill $\square $\medskip

\noindent \textbf{4.12.~Theorem.}~\textit{The following hold for any model }$%
\mathcal{M}$ \textit{of }$\mathrm{PA}$\textit{\ of any cardinality}:$%
\medskip $

\noindent \textbf{(a)} $\mathcal{M}$ \textit{has an expansion to }$\mathrm{CT%
}(c)$\textit{\ iff }$\mathcal{M}$ \textit{has an expansion to} $\mathrm{PAI}%
+\alpha \mathrm{.}$\ \medskip

\noindent \textbf{(b)} $\mathcal{M}$ \textit{has an expansion to }$\mathrm{CT%
}$\textit{\ iff }$\mathcal{M}$ \textit{has an expansion to} $\mathrm{PAI}%
^{\circ }+\alpha \mathrm{.}$\ \medskip

\noindent \textbf{Proof.~}We only verify (a) since the argument for (b) is
similar. The\textbf{\ }left-to-right direction of (a) is justified by
Theorem 4.11. The other direction is evident thanks to the axiom $\alpha .$%
\hfill $\square $\medskip

\noindent \textbf{4.13.~Corollary.~}\textit{The arithmetical consequences of
}$\mathrm{PAI}^{\circ }+\alpha $ \textit{and} $\mathrm{PAI}+\alpha $ \textit{%
are axiomatized by }$\mathrm{PA}+\mathrm{RFN}^{\varepsilon _{0}}(\mathrm{PA}%
).$\medskip

\noindent \textbf{Proof.~}This follows from putting the completeness theorem
of first order logic together with Theorems 4.9 and 4.12.\hfill $\square $%
\medskip

\noindent \textbf{4.14}.~\textbf{Remark.}~In contrast to Theorem 4.6, in
Theorem 4.12 $\mathrm{CT}(c)$ cannot be weakened to $\mathrm{CT}_{0}(c)$,
where $\mathrm{CT}_{0}(c)$ is the fragment of $\mathrm{CT}(c)$ in which the
extended induction scheme is limited to $\Delta _{0}(T)$-formulae. This is
because the arithmetical consequences of $\mathrm{CT}_{0}$ form a tiny
fragment of the arithmetical consequences of $\mathrm{CT}$. More explicitly,
it has been long known by the cognoscenti that by combining results of
Kotlarski \cite{Heniek's CT0} and Smory\'{n}ski \cite{Smorynski (omega-con)}%
, the arithmetical consequences of $\mathrm{CT}_{0}$ can be shown to be
axiomatized by $\mathrm{RFN}^{\omega }(\mathrm{PA})$. The recent work of \L e%
\l yk \cite{Mateusz (reflection)} provides a model-theoretic proof of this
axiomatizability result, and culminates earlier results obtained in
Kotlarski's aforementioned paper, Wcis\l o and \L e\l yk's \cite%
{Bartek+Mateusz(RSL)}, and \L e\l yk's doctoral dissertation \cite{Mateusz
thesis}. \medskip

\noindent \textbf{4.15.~Corollary.~}\textit{The following are equivalent for
every countable model }$\mathcal{M}$\textit{\ of }$\mathit{\mathrm{PA}}$%
:\medskip

\noindent $(i)$ $\mathcal{M}$\textit{\ has an expansion to a model of }$%
\mathrm{PAI}+\alpha .$\medskip

\noindent $(ii)$ $\mathcal{M}$\textit{\ is recursively saturated and
satisfies }$\mathrm{RFN}^{\varepsilon _{0}}(\mathrm{PA}).$\medskip

\noindent \textbf{Proof.~}This is an immediate consequence of Theorems 4.9
and 4.12, together with the resplendence property of countable recursively
saturated models.\hfill $\square $\medskip

\noindent \textbf{4.16}.~\textbf{Remark.}~The proof of Corollary 4.15 makes
it clear that countability restriction on $\mathcal{M}$ can be dropped for $%
(i)\Rightarrow (ii).$ However, the existence of Kaufmann models together
with Tarski's undefinability of truth theorem makes it evident that the
countability restriction cannot be dropped for $(iii)\Rightarrow (i).$
\bigskip

\begin{center}
\textbf{5. ~INTERPRETABILITY\ ANALYSIS OF }$\mathrm{PAI}$\textbf{\ \bigskip }
\end{center}

In this section we examine $\mathrm{PAI}$ through the lens of
interpretability theory, a lens that brings both the semantic and syntactic
features of the theories under its scope into a finer focus. \medskip

\noindent \textbf{5.1.~Theorem.} (mutual $\omega $-interpretability
results). \medskip

\noindent \textbf{(a)} $\mathrm{PAI}$, $\mathrm{UTB}_{0}(c),$ \textit{and} $%
\mathrm{UTB}(c)$ \textit{are pairwise mutually} $\omega $-\textit{%
interpretable}.\medskip

\noindent \textbf{(b)} $\mathrm{PAI}+\alpha $ \textit{and} $\mathrm{CT}(c)$
\textit{are pairwise mutually} $\omega $-\textit{interpretable}.\medskip

\noindent \textbf{(c)} $\mathrm{PAI}^{\circ }+\alpha $ \textit{and} $\mathrm{%
CT}$ \textit{are pairwise mutually} $\omega $-\textit{interpretable}.\medskip

\noindent \textbf{Proof.}~(a) is easy to show using the proofs of Theorem
4.1 and the $(ii)\Rightarrow (iii)$ direction of Theorem 4.6 (note that the
`infinite constant' $c$ in $\mathrm{UTB}_{0}(c),$ and $\mathrm{UTB}(c)$ can
be readily defined in a model of $\mathrm{PAI}$ as \textquotedblleft the
least element of $I$\textquotedblright ). Note that (b)\ follows from
Theorem 4.14. The proof of (c) is similar to the proof of (b).\hfill $%
\square $\medskip

We need some general definitions and results before presenting other results
about $\mathrm{PAI}$. The following definition is motivated by the work of
Albert Visser \cite{Albert-Tehran}; it was introduced in \cite{Ali-Visserian}%
.\medskip

\noindent \textbf{5.2.~Definition.~}Suppose $U$ is a first order theory.
\medskip

\noindent \textbf{(a)} $U$ is \textit{solid} iff the following property $(%
\mathbb{S})$ holds for all models $\mathcal{M}$, $\mathcal{M}^{\ast },$ and $%
\mathcal{N}$ of $U$:\medskip

\noindent $(\mathbb{S})~~~$If $\mathcal{M}\trianglerighteq \mathcal{N}%
\trianglerighteq \mathcal{M}^{\ast }$ and there is an $\mathcal{M}$%
-definable isomorphism\textit{\ }$i_{0}:\mathcal{M}\rightarrow \mathcal{M}%
^{\ast }$, then there is an $\mathcal{M}$-definable isomorphism\textit{\ }$i:%
\mathcal{M}\rightarrow \mathcal{N}$.\medskip

\noindent \textbf{(b)} $U$ is \textit{nowhere solid} if $(\mathbb{S})$\ is
false at every model $\mathcal{M}$ of $U$, i.e., for every model $\mathcal{M}
$ of $U$ there exist models $\mathcal{M}^{\ast },$ and $\mathcal{N}$ of $U$
such that $\mathcal{M}\trianglerighteq \mathcal{N}\trianglerighteq \mathcal{M%
}^{\ast }$ and there is an $\mathcal{M}$-definable isomorphism\textit{\ }$%
i_{0}:\mathcal{M}\rightarrow \mathcal{M}^{\ast }$, but there is no $\mathcal{%
M}$-definable isomorphism\textit{\ }$i:\mathcal{M}\rightarrow \mathcal{N}$%
.\medskip

\noindent Visser showed that \textrm{PA} is a solid theory. The following
proposition can be readily established using the definitions involved (the
proof is straightforward, but notationally complicated).\medskip

\noindent \textbf{5.3.~Proposition.~}\textit{If} $U$ \textit{is a solid
theory and} $V$ \textit{is a retract of} $U$, \textit{then} $V$ \textit{is
also solid. In particular, solidity is preserved by bi-interpretations.}%
\medskip

\noindent The proof of solidity of \textrm{PA} shows the more general
theorem below that will be useful (the proof of Theorem 5.4 is a slight
variant of the proof of solidity of PA presented in \cite{Ali-Visserian}; we
present it for the convenience of the reader).\medskip

\noindent \textbf{5.4.~Theorem.~}\textit{Suppose} $\mathcal{M}_{1}$, $%
\mathcal{M}_{2}$, \textit{and} $\mathcal{M}_{3}$ \textit{are models of} $%
\mathrm{PA}$, \textit{and that }$\mathcal{M}_{i}^{+}$\textit{\ is an} $%
\mathcal{L}_{i}$-\textit{structure that is an expansion of }$\mathcal{M}_{i}$
\textit{and }$\mathcal{M}_{i}^{+}\models \mathrm{PA}(\mathcal{L}_{i})$ (%
\textit{for} $i=1,2,3$). \textit{Then} $(\mathbb{S}^{+})$ \textit{below holds%
}:\medskip

\noindent $(\mathbb{S}^{+})~~~$\textit{If }$\mathcal{M}_{1}^{+}%
\trianglerighteq \mathcal{M}_{2}^{+}\trianglerighteq \mathcal{M}_{3}^{+}$
\textit{and there is an }$\mathcal{M}_{1}^{+}$-\textit{definable
isomorphism\ }$i_{0}:\mathcal{M}_{1}\rightarrow \mathcal{M}_{3}$, \textit{%
then there is an }$\mathcal{M}_{1}^{+}$-\textit{definable isomorphism\ }$i:%
\mathcal{M}_{1}\rightarrow \mathcal{M}_{2}$.\footnote{%
Note that the conclusion \textit{\ }$i:\mathcal{M}_{1}\rightarrow \mathcal{M}%
_{2}$ cannot in general be strengthened to \textit{\ }$i:\mathcal{M}%
_{1}^{+}\rightarrow \mathcal{M}_{2}^{+}$, e.g., let $\mathcal{M}_{1}^{+}=(%
\mathcal{M},D_{1})$ and $\mathcal{M}_{1}^{+}=(\mathcal{M},D_{2})$, where $%
D_{1}$ and $D_{2}$ are distinct $\mathcal{M}$-definable subset of $M$.}
\medskip

\noindent \textit{Consequently,} \textit{an isomorphic copy of} $\mathcal{M}%
_{2}^{+}$ is $\omega $-\textit{interpretable in }$\mathcal{M}_{1}^{+}$ (%
\textit{moreover, the isomorphism at work is} $\mathcal{M}_{1}^{+}$-\textit{%
definable}).\medskip

\noindent \textbf{Proof.~}Suppose $\mathcal{M}_{1}^{+}$, $\mathcal{M}%
_{2}^{+},$ and $\mathcal{M}_{3}^{+}$ are as in the assumption of the
theorem. Further, assume that:

\begin{center}
$\mathcal{M}_{1}^{+}\trianglerighteq \mathcal{M}_{2}^{+}\trianglerighteq
\mathcal{M}_{3}^{+}$,
\end{center}

\noindent and suppose there is an $\mathcal{M}_{1}^{+}$-definable isomorphism%
\textit{\ }$i_{0}:\mathcal{M}_{1}\rightarrow \mathcal{M}_{3}.$ A key
property of $\mathrm{PA}(\mathcal{L})$ is that if $\mathcal{M}^{+}\models
\mathrm{PA}(\mathcal{L})$ and $\mathcal{N}$ is a model of the fragment $%
\mathrm{Q}$ (Robinson's arithmetic) of $\mathrm{PA}$, then as soon as $%
\mathcal{M}^{+}\trianglerighteq \mathcal{N},$ there is an $\mathcal{M}^{+}$%
-definable initial embedding \textit{\ }$j:\mathcal{M}\rightarrow \mathcal{N}
$, i.e., an embedding $j$ such that the image $j(\mathcal{M})$\ of $\mathcal{%
M}$ is an initial submodel of $\mathcal{N}$ (where $\mathcal{M}$ is the $%
\mathcal{L}_{\mathrm{A}}$-reduct of $\mathcal{M}^{+}).$ Hence there is an $%
\mathcal{M}_{1}^{+}$-definable initial embedding \textit{\ }$j_{0}:\mathcal{M%
}_{1}\rightarrow \mathcal{M}_{3}$ and an\textit{\ }$\mathcal{M}_{2}^{+}$%
-definable initial embedding \textit{\ }$j_{1}:\mathcal{M}_{1}\rightarrow
\mathcal{M}_{3}$.\medskip

We claim that both $j_{0}$ and $j_{1}$ are surjective. To see this, suppose
not. Then $j(\mathcal{M}_{1})$ is a proper initial segment of $\mathcal{M}%
_{3}$, where $j$ is the $\mathcal{M}_{1}^{+}$-definable embedding \textit{\ }%
$j:\mathcal{M}_{1}\rightarrow \mathcal{M}_{3}$ given by $j:=j_{1}\circ
j_{0}. $ But then $i_{0}^{-1}(j(M_{1}))$ is a proper $\mathcal{M}_{1}^{+}$%
-definable initial segment of $\mathcal{M}$ with no last element. This is a
contradiction since $\mathcal{M}_{1}^{+}$ is a model of $\mathrm{PA}(%
\mathcal{L}_{1})$, and therefore no proper initial segment of $M_{1}$ is $%
\mathcal{M}_{1}^{+}$-definable. Hence $j_{0}$ and $j_{1}$ are both
surjective; in particular $j_{0}$ serves as the desired $\mathcal{M}_{1}^{+}$%
-definable isomorphism between $\mathcal{M}_{1}$ and $\mathcal{M}_{2}$.
\medskip

Since by assumption $\mathcal{M}_{1}^{+}\trianglerighteq \mathcal{M}%
_{2}^{+}, $ the $\mathcal{M}_{1}^{+}$-definable isomorphism $j_{0}$ allows
us to construct, definably in $\mathcal{M}_{1}^{+}$, an isomorphic copy of $%
\mathcal{M}_{2}^{+}$ whose $\mathcal{L}_{\mathrm{A}}$-reduct is $\mathcal{M}%
_{1}$. Consequently, an isomorphic copy of $\mathcal{M}_{2}^{+}$ is $\omega $%
-interpretable in\textit{\ }$\mathcal{M}_{1}^{+}$. More specifically, the $%
\omega $-interpretation $\mathcal{I}$ at work has the same universe and
arithmetical operations as $\mathcal{M}_{1}$, and for each $n$-ary relation
symbol $R\in \mathcal{L}_{2}\backslash \mathcal{L}_{1,}$ the $\mathcal{I}$%
-interpretation $R^{\mathcal{I}}$ of $R$ is given by $(x_{1},\cdot \cdot
\cdot ,x_{n})\in R^{\mathcal{I}}$ iff $\mathcal{M}_{2}^{+}\models
R(j_{0}(x_{1}),\cdot \cdot \cdot ,j_{0}(x_{n}))$ (the $\mathcal{I}$%
-interpretation of function symbols is analogously defined).\hfill $\square $%
\medskip

\noindent \textbf{5.5.~Theorem.}~\textit{If} $U\in \{\mathrm{PAI}^{\circ },$
$\mathrm{PAI}$, $\mathrm{PAI}+\alpha \},$ \textit{then} $U$\ \textit{is}
\textit{nowhere solid} (\textit{a fortiori}: $U$ \textit{is not solid}%
).\medskip

\noindent \textbf{Proof.~}We present the argument for $U=$ $\mathrm{PAI}$,
similar arguments work for $\mathrm{PAI}^{\circ }$, and $\mathrm{PAI}+\alpha
$. Suppose $(\mathcal{M},I)\models \mathrm{PAI}$, and let $2I=\{2i:i\in I\}.$
The diagonal indiscernibility property of $I$ makes it evident that $(%
\mathcal{M},2I)$ is also a model of $\mathrm{PAI}.$ Note that:

\begin{center}
$(\mathcal{M},I)\trianglerighteq (\mathcal{M},2I)\trianglerighteq (\mathcal{M%
},I)$.
\end{center}

\noindent Clearly the identity map serves as an $(\mathcal{M},I_{1})$%
-definable isomorphism\textit{\ }$i_{0}:(\mathcal{M},I_{1})\rightarrow (%
\mathcal{M},I_{1})$. However, there is no $(\mathcal{M},I)$-definable
isomorphism $f$ between $(\mathcal{M},I)$ and $(\mathcal{M},2I)$ since any
model of $\mathrm{PA(}I\mathrm{)}$ is definably rigid and thus any such
purported $f$ is the identity function, which it impossible for $f$ to be an
isomorphism between $(\mathcal{M},I)$ and $(\mathcal{M},2I)$. This shows
that $\mathrm{PAI}$ is not solid.\hfill $\square $\medskip

\noindent \textbf{5.6.~Theorem.~}\textit{If} $U\in \{\mathrm{PAI}^{\circ },$
$\mathrm{PAI}$ \textit{and} $\mathrm{PAI}+\alpha \},$ \textit{then} $\mathrm{%
PA}$\textit{\ and }$U$\textit{\ are not retracts of each other }(\textit{a
fortiori: }$\mathrm{PA}$ \textit{and }$U$\textit{\ are not bi-interpretable})%
\textit{. }\medskip

\noindent \textbf{Proof.}~Again, we present the argument for $U=$ $\mathrm{%
PAI}$, similar arguments work for $\mathrm{PAI}^{\circ }$, and $\mathrm{PAI}%
+\alpha $. Proposition 5.4 together with Theorem 5.5 show that $\mathrm{PAI}$
is not a retract of $\mathrm{PA}$. To see that $\mathrm{PA}$\ is not a
retract of $\mathrm{PAI}$, it suffices to show that if some model $\mathcal{M%
}$ of $\mathrm{PA}$ is a retract of a model of $\mathrm{PAI}$, then by
Theorem 5.4, $\mathcal{M}$ can parametrically define a class $I$ of
indiscernibles for itself. But this contradicts Corollary 4.3(a).\hfill $%
\square $\medskip

\noindent \textbf{5.7.~Theorem.~}$\mathrm{CT}$ \textit{is solid}. \medskip

\noindent \textbf{Proof.} The solidity of $\mathrm{CT}$ can be established
with the help of Theorem 5.4 and the well-known fact, that if $(\mathcal{M}%
,T_{1},T_{2})\models \mathrm{PA}(T_{1},T_{2})$, where $(\mathcal{M},T_{1})$
and $(\mathcal{M},T_{2})$ are both models of $\mathrm{CT}$, then $%
T_{1}=T_{2} $ (the proof is based on a simple induction, taking advantage of
the assumption that both $T_{1}$ and $T_{2}$ satisfy Tarski's recursive
clauses for all arithmetical formulae in $\mathcal{M}$). \hfill $\square $%
\medskip

\noindent \textbf{5.8.~Remark}.\textbf{~}It is not hard to see that the none
of the theories $\mathrm{CT}(c)\mathrm{,}$ $\mathrm{UTB}$, and $\mathrm{UTB}%
(c)$ are solid. However, $\mathrm{CT}(c)$ has consistent solid extensions.
For example, consider the extension of $\mathrm{CT}(c)$ given by $\mathrm{CT}%
(c)$ + \textquotedblleft $\mathrm{CT}$ is inconsistent\textquotedblright\ +
\textquotedblleft $c$ is the length of the shortest proof of inconsistency
of $\mathrm{CT}$\textquotedblright . By G\"{o}del's second incompleteness
this theory is consistent, and by a reasoning very similar to the proof of
Theorem 5.7 it is also solid.\medskip

\noindent \textbf{5.9.~Theorem.~}$\mathrm{CT}$\ \textit{is a retract of }$%
\mathrm{PAI}^{\circ }+\alpha $, \textit{and} $\mathrm{CT}(c)$\ \textit{is a
retract of }$\mathrm{PAI}+\alpha .$ \medskip

\noindent \textbf{Proof.}~By Theorem 4.14, there is a (uniform) $\omega $%
-interpretation $\mathcal{I}_{c}$ of a model $(\mathcal{M},I^{>c})$ of $%
\mathrm{PAI}+\alpha $ within any model $(\mathcal{M},T,c)$ of $\mathrm{CT}%
(c).$ On the other hand, the definition of $\mathrm{PAI}+\alpha $ makes it
clear that there is a (uniform) $\omega $-interpretation $\mathcal{J}$ of $(%
\mathcal{M},T)$ within $(\mathcal{M},I^{>c})$. A slight variant of this
argument (without the use of the constant $c$) shows that $\mathrm{CT}$\ is
a retract of\textit{\ }$\mathrm{PAI}^{\circ }+\alpha ,$ but we need a
variation of the interpretation $\mathcal{I}_{c}$ in order to show that $%
\mathrm{CT}(c)$ is a retract of $\mathrm{PAI}+\alpha $ because it not clear
that the element $c$ is definable in $(M,I^{>c})$. We can get around this
problem by modifying the interpretation $\mathcal{I}$ as follows: Given any
model $(\mathcal{M},T,c)$ of $\mathrm{CT}(c),$ we first define $I^{>c}$ as
in the interpretation $\mathcal{I}$, and then we define the modified
interpretation $\mathcal{I}_{c}$ given by:

\begin{center}
$\mathcal{I}_{c}(\mathcal{M},T,c)=(\mathcal{M},I_{c},J),$ where $%
I_{c}=\{\left\langle c,i\right\rangle :i\in I^{>c}\}.$
\end{center}

\noindent Thanks to the diagonal indiscernibility property of $I^{>c}$, $(%
\mathcal{M},I_{c})$\ is a model of $\mathrm{PAI}$. Moreover, $(\mathcal{M}%
,I_{c})$ can be shown to be a model of $\mathrm{PAI}+\alpha $ with the same
argument used in the proof of Theorem 4.14 (relying on Lemma 4.11). We can
now readily define an interpretation $\mathcal{J}^{\prime }$ that inverts $%
\mathcal{I}_{c}$ by letting $\mathcal{J}^{\prime }(\mathcal{M},I_{c})=(%
\mathcal{M},T,c)$, where $T$ is the unique truth predicate corresponding to
the partial satisfaction class given by the formula $\sigma $ (of Theorem
4.1), and $c$ is defined as \textquotedblleft the first coordinate of the
ordered pair canonically coded by any member of $I_{c}$\textquotedblright .
Thus $\mathrm{CT}(c)$\ is a retract of\textit{\ }$\mathrm{PAI}+\alpha .$%
\medskip

\noindent \textbf{5.10.~Theorem.~}$\mathrm{PAI}^{\circ }+\alpha $ \textit{is
not a retract of} $\mathrm{CT.}$\medskip

\noindent \textbf{Proof.}~We begin with observing that Proposition 5.3 and
Theorem 5.5 show that the $\omega $-interpretation $\mathcal{I}_{\sigma }$
of $\mathrm{CT}$ in $\mathrm{PAI}^{\circ }$ given by the formula $\sigma $
of Theorem 4.1 is not `invertible', in the sense that there is no
interpretation of $\mathrm{PAI}$ in $\mathrm{CT}$ such that $\mathcal{I}%
_{\sigma }$ and $\mathcal{J}$\ witness that $\mathrm{PAI}$ is retract of $%
\mathrm{CT}$. We next note that if there are interpretations $\mathcal{I}$
and $\mathcal{J}$ that witness that $\mathrm{PAI}+\alpha $ is a retract of $%
\mathrm{CT}$\textrm{, }then Theorem 5.4 assures us that verifiably in $%
\mathrm{PAI}^{\circ },$ the interpretation $\mathcal{I}$ is the same as $%
\mathcal{I}_{\sigma }$ up to a definable permutation of the ambient
universe. This shows that $\mathcal{I}$ is not invertible either, thus
concluding the proof. \hfill $\square $\medskip

\noindent \textbf{5.11.~Question.~}\textit{Is}\textbf{\ }$\mathrm{PAI}%
+\alpha $ \textit{is a retract of} $\mathrm{CT}(c)$?\medskip

\noindent \textbf{5.12.~Question.~}\textit{Is either of the pair of theories
}$\{\mathrm{PAI},$ $\mathrm{UTB}(c)\}$ \textit{a retract of the other one}%
?\medskip

\noindent \textbf{5.13.~Question.~}\textit{Is either of the pair of theories
}$\{\mathrm{PAI}+\alpha ,$ $\mathrm{CT}(c)\}$ \textit{a retract of the other
one}?\medskip

\begin{itemize}
\item We have not succeeded in ruling out that $\mathrm{PAI}$ and $\mathrm{%
UTB}(c)$ are not bi-interpretable (ditto for $\mathrm{PAI}+\alpha $ and $%
\mathrm{CT}(c)$)$.$ We conjecture that the above questions all have negative
answers.\ As partial evidence for our conjecture, let us observe that
Theorem 5.2 and 5.5 show that the $\omega $-interpretation $\mathcal{I}%
_{\sigma }$ of $\mathrm{UTB}(c)$ in $\mathrm{PAI}$, and $\mathrm{CT}(c)$ in $%
\mathrm{PAI}+\alpha ,$ given by the formula $\sigma $ of Theorem 4.1 is not
`invertible', in the sense that there is no interpretation of $\mathrm{PAI}$
in $\mathrm{UTB}(c)$ such that $\mathcal{I}_{\sigma }$ and $\mathcal{J}$\
witness that $\mathrm{PAI}$ is retract of $\mathrm{UTB}(c)$. $\bigskip $
\end{itemize}

\begin{center}
\textbf{5.~FRAGMENTS OF }$\mathrm{PAI}\bigskip $
\end{center}

In this section we briefly examine the model-theoretic behavior of
subsystems $\mathrm{PAI}_{n}$ ($n\in \omega )$ and $\mathrm{PAI}^{-}$ of $%
\mathrm{PAI}$. The section is concluded with two open questions.\medskip

\noindent \textbf{6.1.~Definition.~}For $n\in \omega $, $\mathrm{PAI}_{n}$
is the subsystem of $\mathrm{PAI}$ in which the extended induction scheme
involving $I$ is weakened to $\Sigma _{n}(I)$-formulae, i.e., the axioms of $%
\mathrm{PAI}_{n}$ consist of $\mathrm{PA}$ plus the fragment $\mathrm{I}%
\Sigma _{n}(I)$ of $\mathrm{PA}(I)$, plus axioms (2) and (3) of Definition
3.1 asserting the unboundedness and indiscernibility of $I$. $\mathrm{PAI}%
^{-}$ is the subsystem of $\mathrm{PAI}_{0}$ with no extended induction
scheme involving $I$, so the axioms of $\mathrm{PAI}^{-}$ consist of $%
\mathrm{PA}$ plus axioms (2) and (3) of Definition 3.1.\medskip

\begin{itemize}
\item Given $\mathcal{M}\models \mathrm{PA}$, it is evident that $(\mathcal{M%
},I)\models \mathrm{PAI}^{-}$ iff $I$ is an unbounded set of indiscernibles
in $\mathcal{M}$; and by Theorem 2.1.2, $(\mathcal{M},I)\models \mathrm{PAI}%
_{0}$ iff $I$ is a piecewise-coded unbounded set of indiscernibles in $%
\mathcal{M}$.\medskip
\end{itemize}

\noindent \textbf{6.2.~Theorem.~}\textit{Every model of }$\mathrm{PA}$
\textit{has an elementary end extension that has} \textit{an expansion to}
\textit{a model of} $\mathrm{PAI}_{0}$, \textit{but not to} \textit{a model
of} $\mathrm{PAI.}$ \medskip

\noindent \textbf{Proof.}~Recall that a type $p(x)$ is said to be an
`unbounded indiscernible type' if it is a nonprincipal type satisfying: (1)
there is no constant Skolem term $c$ such that $x\leq c$ is in $p(x)$, and
(2) for any model $\mathcal{M}$ of $\mathrm{PA}$, if $I\subseteq M$ is a set
of elements each realizing $p(x)$, then $I$ is a set of indiscernibles in $%
\mathcal{M}$. It is well-known that unbounded indiscernible types exist in
abundance (continuum-many), and that a type $p(x)$ is minimal (in the sense
of Gaifman) iff $p(x)$ is an unbounded indiscernible type (see theorems
3.1.4 and 3.2.10 of \cite{Roman-Jim's book}). Fix a minimal type $p(x)$ and
any model $\mathcal{M}_{0}$ of $\mathrm{PA}$, and let $\mathcal{M}$ be an $%
\omega $-canonical extension of $\mathcal{M}_{0}$ using $p(x)$ as in section
3.3 of \cite{Roman-Jim's book}. Thus $\mathcal{M}$ is obtained by an $\omega
$-iteration of the process of adjoining an element satisfying $p(x)$. Since $%
p(x)$ is an unbounded indiscernible type, this makes it clear that $\mathcal{%
M}$ carries an unbounded indiscernible subset $I$, and additionally the
order-type of $I$ is $\omega .$ The latter feature makes it clear that $I$
is piecewise-coded in $\mathcal{M}$, and thus $(\mathcal{M},I)\models
\mathrm{PAI}_{0}$ in light of Theorem 2.1.2. \medskip

It remains to show that $\mathcal{M}$ does not have an expansion to $\mathrm{%
PAI.}$ Suppose not, and let $(\mathcal{M},I)\models \mathrm{PAI}.$ It is
easy to see that from the point of view of $(\mathcal{M},I)$, the order-type
of $(I,<)$ is the same as the order-type of $(M,<)$ (where $<$\ is the
ordering on $M$ given by $\mathcal{M}$). Recall that $\mathcal{M}$ can be
written as the union of elementary initial submodels $\mathcal{M}_{n}$ of $%
\mathcal{M}$ (as $n$ ranges in $\omega )$, where $\mathcal{M}_{n}$ is
obtained by $n$-repetitions of the process of adjoining an element
satisfying $p(x)$. By minimality of $p(x)$ this assures us that:\medskip

\noindent $(\ast )$ For each $n\in \omega $, and each choice of $c_{n}\in
M_{n+1}\backslash M_{n}$, $(\mathcal{M}_{n},c_{1},\cdot \cdot \cdot ,c_{n})$
is pointwise definable$.$ \medskip

\noindent The existence of the above isomorphism $f$ makes it clear that
there is some $k\geq 1$ such that $I\cap \left( M_{k}\backslash
M_{k-1}\right) $ is infinite. In particular we can pick distinct $i_{1}$ and
$i_{2}$ in $I\cap \left( M_{k}\backslash M_{k-1}\right) $, together with
elements $\left\{ c_{s}:1\leq s\leq k\right\} $ such that $c_{s}\in
M_{s+1}\backslash M_{s}$ for each $s,$ and moreover $c_{k}$ is below both $%
i_{1}$ and $i_{2}$ (since $M_{k}\backslash M_{k-1}$ has no least element).
By the diagonal indiscernibility property of $I$, this implies that $i_{1}$
and $i_{2}$ are indiscernible in the structure $(\mathcal{M}_{k},c_{1},\cdot
\cdot \cdot ,c_{k},m)_{m\in M_{0}}.$ But this indiscernibility contradicts $%
(\ast )$, and thereby completes the proof.\hfill $\square $\medskip

\noindent \textbf{6.3.~Theorem.~}\textit{If} $\mathcal{M}$ \textit{is a
model of countable cofinality of }$\mathrm{PA}$ \textit{that is expandable
to a model of }$\mathrm{PAI}^{-}$, \textit{then} $\mathcal{M}$ \textit{is
expandable to} \textit{a model of} $\mathrm{PAI}_{0}.$ \textit{However,
every countable model of }$\mathrm{PA}$\textit{\ has an uncountable
elementary end extension that is expandable to a model of }$\mathrm{PAI}^{-}$%
, \textit{but not to }$\mathrm{PAI}_{0}.$\medskip

\noindent \textbf{Proof.}~Suppose $(\mathcal{M},I)\models \mathrm{PAI}^{-},$
where $\mathcal{M}$ has countable cofinality. The countable cofinality of $%
\mathcal{M}$ allows us to construct an unbounded subset $I_{0}$ of $I$ of
order type $\omega $. Since every subset of $\mathcal{M}$ of order-type $%
\omega $ is piecewise-coded, by Theorem 2.1.2, $(\mathcal{M},I_{0})\models
\mathrm{PAI}_{0}.$ To demonstrate the second assertion of the theorem, let $%
\mathcal{M}_{0}$ be a countable model of $\mathrm{PA}$ and $\mathcal{M}$ be
an $\omega _{1}$-canonical extension of $\mathcal{M}_{0}$ using some minimal
type $p(x),$ i.e., $\mathcal{M}$ is obtained by an $\omega _{1}$-iteration
of the process of adjoining an element satisfying $p(x)$ (as in \cite[%
Section 3.3]{Roman-Jim's book}). By Theorem 2.2.14 of \cite{Roman-Jim's book}%
, $\mathcal{M}$ is rather classless, i.e., every piecewise-coded subset of $%
M $ is definable in $\mathcal{M}$. Thus if $(\mathcal{M},I)\models \mathrm{%
PAI}_{0}$, then $(\mathcal{M},I)\models \mathrm{PAI}$, and $I$ is $\mathcal{M%
}$-definable, which contradicts Corollary 4.3(a).\hfill $\square \medskip $

\noindent We close the paper with the following open questions:$\medskip $

\noindent \textbf{6.4.~Question.~}\textit{Does Theorem }6.2 \textit{lend
itself to a hierarchical generalization}? \textit{In other words, is it true
that for every }$n\in \omega $, \textit{every model of }$\mathrm{PA}$
\textit{has an elementary end extension that has} \textit{an expansion to}
\textit{a model of} $\mathrm{PAI}_{n}$, \textit{but not to} \textit{a model
of} $\mathrm{PAI}_{n+1}$? (\textit{It is not even clear how to build a model}
$(\mathcal{M},I)$ of $\mathrm{PAI}_{n}$ \textit{for }$n\in \omega $ \textit{%
that is not a model of} $\mathrm{PAI}_{n+1}$.) $\medskip $

\noindent \textbf{6.5.~Question.~}\textit{Is there a model }$\mathcal{M}$%
\textit{\ of }$\mathrm{PA}$\textit{\ such that }$\mathcal{M}$\textit{\ has}
\textit{an expansion to} \textit{a model of} $\mathrm{PAI}_{n}$ \textit{for
each }$n\in \omega $, \textit{but }$\mathcal{M}$\textit{\ has no expansion
to a model of }$\mathrm{PAI}$?$\bigskip $

\noindent \textsc{Department of Philosophy, Linguistics, and the Theory of
Science \newline
\noindent University of Gothenburg, Gothenburg, Sweden}\newline
\noindent \texttt{email: ali.enayat@gu.se}

\end{document}